\documentstyle{amsppt}
\TagsOnRight
\loadbold
\mathsurround=1.2pt
\define\<{\snug}
\redefine\.{\,.}
\newdimen\theight
\redefine\>{\hskip1pt}
\redefine\le{\leqslant}
\redefine\ge{\geqslant}
\define\){)\hskip1pt}
\define\({\hskip1pt(}
\define\op{\operatorname}
\define\pa{\partial}
\define\paw{{{\partial\over \partial w}}}
\define\la{\lambda}
\define\om{\omega}
\define\Om{\Omega}
\define\a{\alpha}
\define\ga{\gamma}
\define\wt{\widetilde}
\define\U{\text{\rm U}}
\define\bcdot{\,\boldsymbol\cdot\,}
\define\bdot{\boldsymbol\cdot}
\define\opl{\operatornamewithlimits}
\define\itm#1#2{\smallskip\noindent{\rlap{\rom{(#1)}}\qquad}#2\smallskip}
\define\Y#1{{\kern1pt{:}\kern0.5pt#1\kern-0.3pt{:}\kern0.4pt}}
\define\sbh#1{\smallskip\noindent#1\quad}
\define\iso{{\,{\buildrel{{\sim}}\over\longrightarrow}\,}}
\magnification=1100

\topmatter
\title Screenings and a universal Lie-de Rham cocycle
\endtitle
\author Victor Ginzburg and Vadim Schechtman
\endauthor
\address V.G.: Department of Mathematics, University of Chicago, Chicago IL 
60637, USA
\endaddress
\email ginzburg\@math.uchicago.edu
\endemail
\address 
V.S.: Max-Planck-Institut f\"ur Mathematik, Gottfried-Claren-Stra\ss e 26, 
53225 Bonn, Deutschland 
\endaddress
\email vadik\@mpim-bonn.mpg.de\endemail
\dedicatory To A. A. Kirillov on the occasion of his 60th Birthday
\enddedicatory
\endtopmatter

\centerline{\bf Table of Contents}



1.{ $\;$} { Introduction}

{\bf Chapter 1. Toy examples}

2.{ $\;$} { Toy example}

3.{ $\;$} { Generalization}

{\bf Chapter 2. Virasoro algebra}

4.{ $\;$} { Cartan cocycle}

5.{ $\;$} { Bosonization for the Virasoro algebra}

6.{ $\;$} { Bosonic vertex operators}

7.{ $\;$} { Screening operators}

{\bf Chapter 3. $\widehat{\frak{sl}}_2$-case}

8.{ $\;$} { Wakimoto realization}

9.{ $\;$} { The screening current}

{\bf Chapter 4. Affine Lie algebras} 

10.{ $\;$} { Bosonization}

11.{ $\,$} { Screening currents}


\head\S1.\enspace  Introduction\endhead

In the pioneering paper ~\cite{FF}, Feigin and Fuchs  have constructed
intertwining operators between  "Fock-type" modules
over the Virasoro algebra via contour integrals 
of certain operator-valued 
one dimensional local systems over top homology classes
of a configuration space. 
Similar constructions exist for affine Lie algebras. Key
ingredients in such a construction are the so called "screening operators".  
The main observation of the present paper is that the
screening operators contain more information.
 Specifically, at  the chain
level, the screening operators
provide a certain canonical cocycle of the Virasoro (resp. 
affine) Lie algebra with coefficients in the de Rham complex of an
operator-valued local system on the  configuration 
space. This way we obtain   
canonical morphisms from {\it higher} homology 
spaces of the above local systems to  appropriate higher
$Ext$-groups between the  Fock space representations. 

The screening operators that we are interested in this paper are linear
maps $S: M_1 \to M_2[[z,z^{-1}]]$, where $M_1$ and $M_2$ are certain
modules over the Lie algebra ${\frak g}$
in question, e.g., an affine Lie algebra.
We think of such an operator as a "function" of the formal
variable $z$, and write it as $S(z)$. The key property of
screening operators is that they satisfy an equation of the type
$$ [x, S(w)] = {\partial S(x, w)\over \partial w}\quad,\quad
\forall x \in {\frak g}\,,\tag 0.1$$
where $S(x, w):  M_1 \to M_2[[z,z^{-1}]]$ is some other operator,
a companion of $S(w)$. The equation above clearly says that the
pair $S(w),\, S(\bullet, w)$ is a cocycle in the  total complex
associated with a double-complex with Koszul differential, $d'$, and de Rham
differential, $d''$. This observation is a starting point of our
analysis, and the rest of our results is just an elaboration of
that observation. We would like to emphasize that screenings
 still present a  mystery in the sense that we do not
know any kind of general mechanism that would give rise
to operators $S$ satisfying (0.1). In all known cases the screening operators
are constructed by hand, using an explicit description of
the modules $M_1$ and $M_2$ in terms of "creation"
 and "annihilation" operators. Such an explicit description is known
as a {\it bosonisation procedure}, and the bulk of the paper is devoted
to writing out the bosonisation procedure for the modules we need,
since the corresponding formulas are spread over the (mostly physics)
literature (see [FeFr1], [FeFr2] for mathematical results).

Our construction is
motivated by, and in a special case reduces to
the construction of ~\cite{BMP1}, \cite{BMP2}. In fact, as the results 
of {\it loc. cit} and ~\cite{FS} suggest, an explanation 
of our construction should lie in some equivalence of (derived) 
categories of representations of quantum groups, and the corresponding 
affine algebras, sending (contragradient) Verma modules 
to Wakimoto modules.    

On the other hand, we believe that the cocycles that 
we study can be most adequately interpreted as
de Rham cohomology classes of the chiral algebras considered
by Beilinson-Drinfeld \cite{BD}. 

Recently A. Sebbar, ~\cite{S},  was able to obtain a 
"$q$-deformation" of  most of the
constructions of this paper, where the de Rham cohomology 
is deformed to Aomoto-Jackson $q$-de Rham cohomology. It is interesting 
to note that, as was expected, the operators from Lemma 3.3 are deformed 
to Kashiwara operators (see ~\cite{Ka} or ~\cite{L}, n.1.213).    

The main results of this paper have been obtained back in 1991, and 
written up in February 1996. We are deeply indepted to Ed Frenkel 
for the numerous enlightning discussions on bosonisation
during 1990-1991. 
We thank Jim Stasheff 
for useful remarks.

\head Chapter 1. Toy examples.\endhead

\head\S2.\enspace  The toy example\endhead

In this work, everything will be over the base field $\Bbb C$.

\sbh{2.1.}Let $\frak{g}$ be the Lie algebra $\op{sl}(2)$, with
the standard
generators $E$, $F$, $H$.
For $\la\in\Bbb C$, let $M(\la)$ denote
the Verma module over $\frak{g}$ generated by the vacuum vector
$v_{\la}$ subject to the relations $Ev_{\la}=0$, $Hv_{\la}=
\la v_{\la}$. We shall use the formulas
$$
EF^av_{\la}=a(\la-a+1\)F^{a-1}v_{\la},\qquad
HF^av_{\la}=(\la-2a)F^av_{\la}\.
$$
(In the sequel, in our formulas we shall use the agreement $F^bv_{\la}=0$
for $b<0$.)

\TagsOnLeft

Let us pick $\la,\la'\in\Bbb C$.
For an integer $n\ge 0$,  we consider the $\Bbb C$-linear operator
$$
V_n\: M(\la'-1)\to M(\la-1),\qquad
V_n(F^av_{\la'-1})=F^{a+n}v_{\la-1}\.
$$
These operators satisfy the following commutation relations.
$$
\align
[E,V_n](F^av_{\la'-1})&=(-n^2+(\la-2a\)n+a(\la-\la'))
F^{a+n-1}v_{\la-1},\tag"(a)"\\
[H,V_n](F^av_{\la'-1})&=(\la-\la'-2n)F^av_{\la-1},\tag"(b)"\\
[F,V_n]&=0\.\tag"(c)"
\endalign
$$

\sbh{2.2.}Consider the operator-valued formal power series
$$
V(z)=\sum_{n\ge 0}V_nz^{-n-1}dz\:\enspace M(\la'-1)@>>> M(\la-1)[[z^{-1}]]\,
\frac{dz}{z}
$$
(here $z$ is a formal variable). Let us try to find a number
$\a\in\Bbb C$ and an operator
$$
V(E;z)=\sum_{n\ge 0} V_n(E\)z^{-n}\:\enspace M(\la'-1)@>>> M(\la-1)[[z^{-1}]]
$$
such that
$$
[E,V(z)]=(d+\a\,dz/z\)V(E;z)\.\tag"(a)"
$$

The equation (a) is equivalent to the system of equations
$$
[E,V_n]=(-n+\a\)V_n(E)\qquad(n\ge 0)\.\tag"(b)"
$$
So, we have
$$
(-n+\a\)V_n(E)(F^av_{\la'-1})=(-n^2+(\la-2a\)n+a(\la-\la')\)
F^{a+n-1}v_{\la-1},
$$
therefore
$V_n(E)(F^av_{\la'-1})=(n+\beta(a))f^{a+n-1}v_{\la-1}$ for some
function
$\beta(a)$ such that
$$
(-n+\a)(n+\beta(a))=-n^2+(\la-2a\)n+a\(\la-\la'),
$$
that is,
$\beta(a)-\a=-\la+2a$, i.e.,
$\beta(a)=2a+\a-\la$, and $\beta(a\)\a=a\(\la-\la')$
for all $a$.

Suppose that $\la\ne\la'$. Then we must have
$\beta(a)=2a$, $\a=\la$, hence from the second equation
we obtain $\la'=-\la$.

\sbh{2.3.}From now on we suppose that $\la'=-\la$. Thus, we have
$$
[E,V_n]=(-n+\la\)V_n(E)\quad\text{where}\quad
V_n(E)(F^av_{-\la-1})=(n+2a\)F^{a+n-1}v_{\la-1}\,.\tag"(a)"
$$

From~2.1 (b) we obtain
$$
[H,V_n]=(-n+\la\)V_n(H),\tag"(b)"
$$
where $V_n(H)=2V_n$. Finally,
$$
[F,V_n]=0\.\tag"(c)"
$$
Therefore, we come to the following conclusion.

\sbh{2.4.}The operator
$$
V(z)=\sum_{n\ge 0}V_nz^{-n-1}dz\:\enspace M(-\la-1)@>>> M(\la-1)[[z^{-1}]]\,
\frac{dz}{z}
$$
defined by $V_n(F^av_{-\la-1})=F^{a+n}v_{\la-1}$, satisfies
the following relation
$$
[X,V(z)]=(d+\la\,dz/z\)V(X;z)\qquad (X\in\frak{g}),\tag"(a)"
$$
where the operators
$$
V(X;z)=\sum_{n\ge 0}V_n(X\)z^{-n}\:\enspace M(-\la-1)@>>> M(\la-1)[[z]],
$$
linearly depending on $X\in\frak{g}$, are defined by
$$
V_n(E)(F^av_{-\la-1})=(n+2a\)F^{a+n-1}v_{\la-1},\quad
V_n(H)=2V_n,\quad V_n(F)=0\.
$$

\sbh{2.5.}We have
$$
\align
[H,V_n(E)](F^av_{-\la-1})&=-2\(n+2a)(n-\la-1)F^{a+n-1}v_{\la-1},\\
[F,V_n(E)]&=-V_n(H),\\
[E,V_n(H)](F^av_{-\la-1})&=-2\(n+2a)(n-\la)F^{a+n-1}v_{\la-1},\\
[F,V_n(H)]&=0\.
\endalign
$$
It follows that for any $X,Y\in\frak{g}$ and $n\ge 0$, we have
$$
V_n([X,Y])=[X,V_n(Y)]-[Y,V_n(X)]\.\tag"(a)"
$$

\sbh{2.6.}Let us
consider the complex (of length $1$)
$$
\Om^{\bdot}\:\enspace 0@>>>\Om^0@>d_{\la}>>\Om^1@>>>0,\tag"(a)"
$$
where $\Om^0=\Bbb C[[z^{-1}]]$, $\Om^1=\Bbb C[[z^{-1}]]\,dz/z$,
$d_{\la}=d+\la\,dz/z$.

We will always write $\op{Hom}$ for $\op{Hom}_{\Bbb C}$
unless specified otherwise. For any two integers $i,j\ge 0$ set
$$
\split
C^{ij}(\frak{g};M(-\la-1),\;M(\la-1))
&=\op{Hom}\(\Lambda^i\frak{g}\otimes M(-\la-1),\;
M(\la-1)\otimes\Om^j)\\
&=\op{Hom}\(\Lambda^i\frak{g},\op{Hom}\(M(-\la-1),M(\la-1)\otimes\Om^j))\.
\endsplit
$$
The bigraded space $C^{\bdot\bdot}(\frak{g}; M(-\la-1),M(\la-1))$
has the natural structure of a bicomplex. The first differential
$$
d'\:\enspace C^{ij}(\frak{g}; M(-\la-1),M(\la-1))\to
C^{i+1,j}(\frak{g}; M(-\la-1),M(\la-1))
$$
is induced by the standard Koszul differential in the cochain complex
of the Lie algebra $\frak{g}$ with coefficients in the module
$\op{Hom}\(M(-\la-1),M(\la-1))$.

The second differential
$$
d''\:\enspace C^{ij}(\frak{g}; M(-\la-1),M(\la-1))\to
C^{i,j+1}(\frak{g}; M(-\la-1),M(\la-1))
$$
is induced by the differential in the complex $\Om^{\bdot}$.

Let $C^{\bdot}(\frak{g};M(-\la-1),M(\la-1))$ denote the associated
total complex.

The operator $V(z)$ is an element of the space $C^{01}(\frak{g};M(-\la-1),
M(\la-1))$. Let us denote this element by $V^{01}(z)$.
The operators $V(X;z)$ (\<$X\in\frak{g}$) define
an element $V^{10}(z)$ of the space $C^{10}(\frak{g};M(-\la-1),M(\la-1))$.

Property 2.4 (a) means that $d'(V^{01}(z))=d''(V^{10}(z))$.
Property 2.5 (a) means that $d''(V^{10}(z))=0$.
Therefore, the element $(V^{01}(z),V^{10}(z))$ is a {\it $1$-cocycle}
of the total complex $C^{\bdot}(\frak{g};M(-\la-1),M(\la-1))$.

\sbh{2.7.}Suppose
that $\la$ is a nonnegative integer.
In this case, the complex
$\Om^{\bdot}$ has two one-dimensional cohomology spaces.
The space
$H^0(\Om^{\bdot})$ is
generated by the function $z^{-\la}$; and the space
$H^1(\Om^{\bdot})$ generated by the class of the form $z^{-\la}\,dz/z$.
(If $\la\not\in\Bbb N$, the complex $\Om^{\bdot}$ is
acyclic.)

Consider the dual spaces $H_i=H^i(\Om^{\bdot})^*$.
The space $H_1$ is generated by the functional $\Om^1\to\Bbb C$
which assigns to a form $\om$ the residue $\opl{res}_{z=0}(\om z^{\la})$.
The space $H_0$ is generated by the (restriction of) the functional
$\Om^0\to\Bbb C$ which assigns to a function $f(z)$ the residue
$\opl{res}_{z=0}(f(z\)z^{\la}\,dz/z)$.

The previous discussion implies the following.

\itm{a}{The operator
$\opl{res}_{z=0}(V^{01}(z\)z^{\la})\in\op{Hom}\(M(-\la-1),
M(\la-1))$ is an intertwiner.
It is the unique $\frak{g}$-homomorphism $M(-\la-1)\to M(\la-1)$ sending
$v_{-\la-1}$ to $F^{\la}v_{\la-1}$.}

\itm{b}{The operator $\opl{res}_{z=0}(V^{10}(z\)z^{\la}\,dz/z)
\in\op{Hom}\(\frak{g},\op{Hom}\(M(-\la-1),M(\la-1)))$
is a $1$-cocycle of the Lie algebra $\frak{g}$
with coefficients in the $\frak{g}$-module
$\op{Hom}\(M(-\la-1),\,M(\la-1))$.}

Therefore, this operator defines certain element of the space
$\op{Ext}^1_{\frak{g}}(M(-\la-1),\,M(\la-1))$.

\head\S3.\enspace  Generalization of the toy example\endhead

\sbh{3.1.} Let $A=(a_{ij})_{i,j=1}^r$ be a symmetrizable generalized 
Cartan matrix,  
and let $\frak{g}$ be the corresponding Kac-Moody Lie algebra defined by the 
Chevalley  
generators $E_i, F_i, H_i$ $(i=1,\ldots, r)$ and relations
(see ~\cite{K}, 0.3)
$$[H_i,H_j]=0\,;\enspace [H_i,E_j]=a_{ij}E_i,\enspace  [H_i,F_i]=-a_{ij}F_i\,;
\enspace  [E_i,F_j]=\delta_{ij}H_i\,;$$

$$\op{ad}\(E_i)^{-a_{ij}+1}(E_j)=\op{ad}\(F_i)^{-a_{ij}+1}(F_j)=0\,.$$

Let $\frak{h}\subset\frak{g}$ be the Cartan subalgebra spanned by the elements 
$H_1,\ldots,H_r$.  
For $i=1,\ldots,r$, let $\alpha_i\in\frak{h}^*$ be the corresponding simple 
root; let $r_i:\frak{h}^*\longrightarrow\frak{h}^*$ be the corresponding 
simple reflection, 
$$
r_i\lambda=\lambda-\langle H_i,\lambda\rangle\alpha_i.
$$
Let $\rho\in\frak{h}^*$ be the element defined by $\langle H_i,\rho\rangle=
1\ (i=1,\ldots,r)$. Let $\frak{n}_-\subset\frak{g}$ be the Lie subalgebra generated by 
the elements $F_1,\ldots,F_r.$

For $\lambda\in\frak{h}^*$, let $M(\lambda)$ denote the Verma module over $\frak{g}$, 
with one generator $v_{\lambda}$ and relations 
$E_iv_{\lambda}=0;\ H_iv_{\lambda}=\langle H_i,\lambda-\rho\rangle 
v_{\lambda}$.

\proclaim{\rom{3.3.}\quad Lemma-definition}
 There exists a unique linear 
operator $\partial_i:\U\frak{n}_-\longrightarrow\U\frak{n}_-$ such that $\partial_i(F_j)=\delta_{ij}\cdot 1
\ (j=1,\ldots,r)$ and for any $x,y\in\U\frak{n}$, $\partial_i(xy)=\partial_i(x)y+x\partial_i
(y)$.\endproclaim

\demo{Proof}
 The uniqueness is clear. Let us prove the existence. Let $A$ 
be the free associative $\Bbb C$-algebra with generators $\theta_j,\ 
j=1,\ldots, r$. It is clear that there exists a unique linear operator 
$\partial_i:A\longrightarrow A$ such that 
$\partial_i(\theta_j)=\delta_{ij}\cdot 1\,,$ and such that,
for any $x,y\in A$, one has $\partial_i(xy)=\partial_i(x)y+x\partial_i(y)$. 

For an integer $a\geq 1$ and $j\neq k$ in $\{1,\ldots,r\}$, define 
the following element in the algebra $A$  
$$
C(j,k;a)=\op{ad}(\theta_j)^a(\theta_k)=\sum_{p=0}^a
(-1)^p\binom ap\theta_j^{a-p}\theta_k\theta_j^p\.
$$
We claim that

\itm{a}{for any $a$, $j$, $k$, $C(j,k;a)\in\op{Ker}\(\pa_i)$.}

Indeed, the claim is clear for $(j,k)$ such that $i\ne j$ and $i\ne k$.
We have
$$
\pa_i(C(j,i;a))=\sum_{p=0}^a(-1)^p\binom ap
\theta_i^a=(1-1)^a\theta_j^a=0\.
$$
Let us prove that $\pa_i(C(i,k;a))=0$ by induction on $a$.
For $a=1$ this is obvious. We have
$$
\align
\pa_i(C(i,k;a))&=\pa_i(\theta_iC(i,k;a-1)-C(i,k;a-1\)\theta_i)\\
\intertext{(by induction)}
&=\pa_i(\theta_i\)C(i,k;a-1)-C(i,k;a-1\)\pa_i(\theta_i)\\
&=C(i,k;a-1)-C(i,k;a-1)=0\.
\endalign
$$
The claim is proved.

It follows from (a) that the operator $\pa_i\:A\to A$ induces an operator
$\pa_i\:\U\frak{n}_-\to\U\frak{n}_-,$
since $\U\frak{n}_-$ is the quotient of $A$
by the two-sided ideal generated by the elements
$C(j,k;-a_{jk}+1)$. The lemma is proven.\qed
\enddemo

\sbh{3.2}{\bf Remark.} The operators $\partial_i$ are classical limits of
the Kashiwara operators, see \cite{Ka},
on the quantized universal enveloping algebra
$\U_q(\frak{n}_-)$. In fact most of the constructions of this and the
subsequent
sections have been recently extended to the quantized setup
in \cite{S}.

\sbh{3.3.} Pick an element  
$\lambda\in\frak{h}^*$ and $i\in\{1,\ldots,r\}$. Set 
$\lambda'=r_i\lambda$. For each integer $n\geq 0$, 
we introduce  linear operators
$V_{i;n}\,,$ and $V_{i;n}(E_j):
\enspace M(\lambda')\longrightarrow M(\lambda)$ defined
for $x\in\U\frak{n}_-$ by 
$$
V_{i;n}:\; xv_{\lambda'} \mapsto xF_i^nv_{\lambda}\quad,\quad
V_{i;n}(E_j):\; xv_{\lambda'} \mapsto 
a_{ji}\partial_j(x)F_i^nv_{\lambda}+
\delta_{ij}nxF_i^{n-1}v_{\lambda}.
$$
We further form
the following elements of $\op{Hom}(M(r_i\lambda),M(\lambda)
[[z^{-1}]]\frac{dz}{z})$:
$$
V_i(z)=\sum_{n=0}^{\infty}V_{i;n}z^{-n-1}dz\quad,\quad
V_i(E_j;z)=\sum_{n=0}^{\infty}V_{i;n}(E_j)z^{-n}\,.
$$

\proclaim{\rom{3.5.}\quad Proposition}  For any $i,j$ and $n\ge 0$, we have
$\enspace 
[E_j,\,V_{i;n}]=$
$(-n+\langle H_i,\lambda\rangle)\cdot$
$V_{i;n}(E_j),\; $
 equivalently, there is a power series identity
$$
[E_j,V_i(z)]=(d+\langle H_i,\lambda\rangle\frac{dz}{z})V_i(E_j;z).
$$\endproclaim

\demo{Proof}
 For any $i$ and $x\in \U\frak{n}_-$
we can write by definition
$$V_i(z) :\enspace  x\cdot v_{\lambda'} \mapsto x e^{z\cdot F_i}\cdot v_\lambda.$$
Let $E_j$ be a simple Chevalley generator.
We will frequently use the following formulas:
$$[E_j, x] = \partial_j(x)\cdot H_j + x_1\quad,\quad  x_1\in \U\frak{n}_-\,. \eqno(a)$$
$$e^{-z\cdot F_i} \cdot H_j \cdot e^{z\cdot F_i}= H_j + \alpha_i(H_j)\cdot z \cdot F_i\eqno(b)$$
$$E_j e^{z\cdot F_i} - e^{z\cdot F_i} E_j =P(z)\in z\cdot\U\frak{n}_-[z],\quad\op{and}\quad
P=0\enspace \op{if}\enspace   i\not= j\eqno(c)$$
Thus, $P$ is a polynomial in $z$ without constant term with values
in  $\U\frak{n}_-$.

Using formulas (a)-(b) above, we find:
$$[E_j, xe^{z\cdot F_i}] = [E_j, x]\cdot e^{z\cdot F_i} + x \cdot
[E_j, e^{z\cdot F_i}]=
(\partial_j(x)\cdot H_j + x_1)\cdot e^{z\cdot F_i} + x \cdot P(z)\,.$$
Hence,
$$E_j\Bigl(V_i(z)(x v_{{\lambda'}})\Bigr)= E_j(xe^{z\cdot F_i}v_{\lambda})
=(\partial_j(x)\cdot H_j + x_1)\cdot e^{z\cdot F_i}v_{\lambda}
x \cdot P(z)v_{\lambda}\tag *$$
$$=\partial_j(x)\cdot e^{z\cdot F_i} (H_j+ z\alpha_i(H_j) F_i)v_{\lambda}
+x_1e^{z\cdot F_i}v_{\lambda} + x \cdot P(z)v_{\lambda}\,.$$
On the other hand, from (a) we get
$$E_jx v_{{\lambda'}}= [E_j, x]\cdot v_{{\lambda'}} 
= \partial_j(x)\cdot H_j\cdot v_{{\lambda'}}
+ x_1
v_{{\lambda'}}\,.$$
Hence,
$$V_i(z)(E_jx v_{{\lambda'}})=\partial_j(x)\cdot
 H_je^{z\cdot F_i}v_{\lambda}
+x_1e^{z\cdot F_i}v_{\lambda}\quad\op{by (b) } \tag **$$
$$=\partial_j(x)\cdot e^{z\cdot F_i}(H_j + 
\alpha_i(H_j)\cdot z \cdot F_i)+x_1e^{z\cdot F_i}v_{\lambda}.$$
From $(*)$ and $(**)$ we obtain
$$[E_j, V_i(z)](x v_{{\lambda'}})=
\langle {\lambda'} -\lambda, H_j\rangle\partial_j(x)\cdot 
e^{z\cdot F_i}v_{\lambda}
+\langle \alpha, H_j\rangle z\cdot F_i e^{z\cdot F_i}v_{\lambda} +
x \cdot P(z)v_{\lambda}\,.$$

It is easy to deduce from the last formula and formula
(c) above, that if $i\not= j$ the RHS of takes
the form  
$$(d+\langle H_i,\lambda\rangle\frac{dz}{z})V_i(E_j;z)$$
The case $i=j$ is treated similarly.\qed
\enddemo

\sbh{3.6.} For each $n\ge 0$ we
 define the operators $V_{i;n}(H_j): M(\lambda')\longrightarrow 
M(\lambda)$ by    $\enspace  V_{i;n}(H_j)=a_{ji}V_{i;n}\enspace  $. 
One checks easily that  

(b) $[H_j,V_{i;n}]=(-n+\langle H_i,\lambda\rangle)V_{i;n}(H_j)$. 

Finally, it is evident that 

(c) $[F_j,V_{i;n}]=0$ for all $j$. 

Set $V_{i;n}(F_j)=0$. 

\proclaim{\rom{3.7.}\quad Lemma-definition}
 For each
$n\ge 0$, there exists a unique element
$$
V_{i;n}(\bcdot)\in\op{Hom}\(\frak{g},\op{Hom}\(M(\la'),\,M(\la)))
$$
such that
$V_{i;n}(E_j)$, $V_{i;n}(H_j)$, $V_{i;n}(F_j)$ are the elements
defined
above and the following
cocycle condition holds\/\rom:

\itm{a}{For any $X,Y\in\frak{g}$,
$V_{i;n}([X,Y])=[X,V_{i;n}(Y)]-[Y,V_{i;n}(X)]$.\qed}
\endproclaim

The previous considerations may be reformulated in terms of the
generating functions
$$
V_i(X;z)=\sum_{n=0}^{\infty}V_{i;n}(X)z^{-n}\in\op{Hom}(M(r_i\lambda),
M(\lambda)
[[z^{-1}]]),\quad X\in\frak{g}.
$$
as the follows. 

\proclaim{\rom{3.9.}\quad Theorem} For any $X,Y\in\frak{g}$,
we have
$$
\align
[X,V_i(z)]&=(d+\langle H_i,\la\rangle\,dz/z\)V_i(X;z),\\
V_i([X,Y];z)&=[X,V_i(Y;z)]-[Y,V_i(X;z)]\.\qed
\endalign
$$
\endproclaim

\sbh{3.10.}
 Suppose an element $w$ of the Weyl group of $\frak{g}$ 
together with its reduced decomposition 
$w=r_{i_a}\cdot\ldots\cdot r_{i_1}$, and an element 
$\lambda\in\frak{h}^*$ is given. 
For $p=1,\ldots,a$, set 
$\lambda_p=r_{i_{p-1}}\cdot\ldots\cdot r_{i_1}\lambda$. 

Define a complex 
$$
\Omega^{\bullet}:\enspace \ 0\longrightarrow\Omega^0\longrightarrow\ldots
\longrightarrow\Omega^a\longrightarrow 0
$$
as follows. Set $A=\Bbb C[[z_1^{-1},\ldots,z_a^{-1}]]$. By definition, 
$\Omega^p$ is the free $A$-module with the basis 
$\{(dz_{j_1}/z_{j_1})\wedge\ldots\wedge (dz_{j_p}/z_{j_p}),\ 1\leq 
j_1<\ldots
< j_p\leq a\}$. The differential is defined by 
$$
d\eta=d_{DR}(\eta)+(\sum_{p=1}^a\frac{\langle H_{i_p},
\lambda_p\rangle}{z_p}
dz_p)\wedge\eta
$$
where $d_{DR}$ is the de Rham differential. 

\sbh{3.11.} For each $p=1,\ldots,a$, consider the operators 
$\omega_p=V_{i_p}(z_p):\enspace \ M(\lambda_{p+1})=M(r_{i_p}\lambda_p)
\to
M(\lambda_p)z_p^{-1}[[z_p^{-1}]]$ and 
$\tau_p(X)=V_{i_p}(X;z_p):\enspace \ M(\lambda_{p+1})\to
M(\lambda_p)[[z_p^{-1}]],$
$X\in\frak{g},$ defined above. 
 
For each $m$, $0\leq m\leq a$, us define the operators 
$$
V^{m,a-m}\in\op{Hom}(\Lambda^m\frak{g},\op{Hom}(M(w\lambda),M(\lambda)
\otimes\Omega^{a-m}))
$$
as follows. We set
 $$
\split
&V^{m,a-m}(X_1\wedge\cdots\wedge X_m)\\
&\qquad=\sum_{1\le p_1<\cdots<p_m\le a}
(-1)^{\op{sgn}(p_1,\dots,p_m)}\\
&\qquad\qquad\times\bigg(\sum_{\sigma\in\Sigma_m}(-1)^{\op{sgn}(\sigma)}
\om_1\cdots\tau_{p_1}(X_{\sigma(1)})\cdots
\tau_{p_m}(X_{\sigma(m)})\cdots\om_a\bigg)\\
&\qquad\qquad\times dz_1\wedge\cdots\wedge\widehat{dz}_{p_1}\wedge\cdots
\wedge\widehat{dz}_{p_m}\wedge\cdots\wedge dz_a\.
\endsplit
$$
Here for each sequence $1\leq p_1<\ldots<p_m\leq a$, the corresponding 
summands are obtained by replacing the operators $\omega_{p_j}$ 
in the product $\omega_1\cdot\ldots\cdot \omega_a$, by 
$\tau_{p_j}(X_{\sigma(j)})$. 
$\Sigma_m$ denotes the group of all bijections $\sigma:\{1,\ldots,m\}
\iso
\{1,\ldots,m\}$.  
The power $\op{sgn}(p_1,\ldots,p_m)$ is defined by induction on $m$ as 
follows. 
We set 
$$
\op{sgn}(\emptyset)=0;\ \op{sgn}(p_1,\ldots,p_m)=\op{sgn}(p_1,\ldots,p_{m-1})+p_m+m.
$$
For example, 
$$
V^{0a}=V_{i_1}(z_1)\cdot\ldots\cdot V_{i_a}(z_a)dz_1\wedge\ldots\wedge
 dz_a.
$$
Consider the double complex $C^{\bullet}(\frak{g};
\op{Hom}(M(w\lambda),
M(\lambda)\otimes\Omega^{\bullet}))$. Here the first differential is 
the Koszul differential in the complex of cochains of the Lie algebra 
$\frak{g}$ with coefficients in the complex of $\frak{g}$-modules
 $\op{Hom}(M(w\lambda),
M(\lambda)\otimes\Omega^{\bullet})$. The action of $\frak{g}$ is
induced by the standard action (by the commutator) of $\frak{g}$ on 
$\op{Hom}(M(w\lambda),M(\lambda))$. The second differential is induced 
by the differential in $\Omega^{\bullet}$. 

By definition, we have 
$$
V^{m,a-m}\in C^m(\frak{g};\op{Hom}(M(w\lambda),M(\lambda)\otimes
\Omega^{a-m})).
$$ 
\proclaim{\rom{3.12.}\quad Theorem} The element $V=(V^{0a},\dots,V^{a0})$
is
an $a$-cocycle in the total complex associated with the double complex
$C^{\bdot}(\frak{g};\op{Hom}\(M(w\la),M(\la)\otimes\Om^{\bdot}))$.\qed
\endproclaim

\proclaim{\rom{3.13.}\quad Corollary} The cocycle $V$ induces linear maps
$$
f_m\:\enspace H^m(\Om^{\bdot})^*\to\op{Ext}^{a-m}_{\frak{g}}(M(w\la),M(\la)),
\qquad m=0,\dots,a\.
$$
\endproclaim

\example{\rom{3.14.}\quad Example}
 Assume that all numbers $\langle H_p,\lambda_{p}
\rangle,\ p=1,\ldots, a,$ are nonnegative integers. The highest homology 
space $H^a(\Omega^{\bullet})$ is one-dimensional, generated by the 
(image of) the functional $r\in\Omega^{a*}$ defined by 
$$
r(\eta)=\op{res}_{z_a=0}\ldots\op{res}_{z_1=0}(z_1^{\langle H_{i_1},
\lambda_1\rangle}
\ldots z_a^{\langle H_{i_a},\lambda_a\rangle}\eta).
$$
The image $f_0(r)\in\op{Hom}_{\frak{g}
}(M(w\lambda),M(\lambda))$ is the unique 
(up to proportionality) intertwiner sending 
$v_{w\lambda}$ to $F_{i_a}^{\langle H_{i_a},\lambda_a\rangle+1}\cdot\ldots
\cdot F_{i_1}^{\langle H_{i_1},\lambda_1\rangle+1}v_{\lambda}$. 
\endexample

\head Chapter 2. Virasoro algebra.\\
\S4.\enspace  Cartan cocycle\endhead

\sbh{4.1.}Let $\frak{g}$ be a Lie algebra. Consider the differential 
graded
Lie algebra $\frak{g}^{\bdot}=\frak{g}^{-1}\oplus\frak{g}^0$
defined as follows. We set
$\frak{g}^{-1}=\frak{g}^0=\frak{g}$;
the differential $d\:\frak{g}^{-1}\to\frak{g}^0$
is the identity map; the bracket $\Lambda^2\frak{g}^0\to\frak{g}^0$
coincides with
the bracket in $\frak{g}$;
the bracket $\frak{g}^{-1}\otimes\frak{g}^0\to\frak{g}^0$ coincides
with the bracket in $\frak{g}$.

For $x\in\frak{g}$,  we denote by the same letter $x$ the corresponding
element
of $\frak{g}^0$, and by $i_x$ the corresponding element of $\frak{g}^{-1}$.

As a dg Lie algebra, the algebra $\frak{g}^{\bdot}$ is generated
by the Lie algebra $\frak{g}$ (in degree $0$) and by
the elements $i_x\in\frak{g}^{-1}$
(\<$x\in\frak{g}$) subject to the relations (a)--(c) below.
$$
\alignat2
d(i_x)&=x&\qquad& (x\in\frak{g}),\tag"(a)"\\
[x,i_y]&=i_{[x,y]}&\qquad&(x,y\in\frak{g}),\tag"(b)"\\
[i_x,i_y]&=0&\qquad& (x,y\in\frak{g})\.\tag"(c)"
\endalignat
$$

\sbh{4.2.}Let $M^{\bdot}$ be a complex of vector spaces.
A $\frak{g}^{\bdot}$-module structure on $M^{\bdot}$ is the same 
thing as a
collection of data (a), (b) below satisfying properties (c)--(e) below.

\itm{a}{Morphisms of complexes $x\:M^{\bdot}\to M^{\bdot}$,
(\<$x\in\frak{g}$)
that define an action of the Lie algebra $\frak{g}$.}

\itm{b}{Morphisms of graded spaces
$i_x\:\enspace M^{\bdot}\to M^{\bdot}[-1]$,
(\<$x\in\frak{g}$).}

\itm{c}{({\it Cartan formula\/}) $[d,i_x]=x$ (\<$x\in\frak{g}$).}

\noindent Here (and below)
the commutators are understood in the graded sense,
i.e., $[d,i_x]=d\circ i_x+i_x\circ d$.

\itm{d}{$[x,i_y]=i_{[x,y]}$ (\<$x,y\in\frak{g}$).}

\itm{e}{$[i_x,i_y]=0$ (\<$x,y\in\frak{g}$).}

\sbh{4.3.}Let us
consider the enveloping algebra $\U\frak{g}^{\bdot}$.
It is a dg associative algebra.
We have the canonical embedding
$$
\Lambda^{\bdot}(\frak{g}^{-1})\hookrightarrow\U\frak{g}^{\bdot}.
$$
We shall use the notation $i_{x_1\dots x_a}$ for the elements
$i_{x_1}\cdots i_{x_a}\in (U\frak{g}^{\bdot})^{-a}$,
\<$x_1,\dots,x_a\in \frak{g}$, where
$(U\frak{g}^{\bdot})^{k}$ stands for grade degree $k$ component
of $U\frak{g}^{\bdot}$.

\proclaim{\rom{4.4.}\quad Lemma}
For any $x_1,\dots,x_a\in\frak{g}$, we have
$$
\split
di_{x_1\dots x_a}&=\sum_{p=1}^a(-1)^{p-1}x_pi_{x_1\dots\hat{x}_p\dots
x_a}+\sum_{1\le p<q\le a}(-1)^{p+q}i_{[x_p,x_q]x_1\dots\hat{x}_p
\dots\hat{x}_q\dots x_a}\\
&=\sum_{p=1}^a(-1)^{p-1}i_{x_1\dots\hat{x_p}\dots x_a}x_p+
\sum_{1\le p<q\le a}(-1)^{p+q+1}i_{[x_p,x_q]x_1\dots\hat{x}_p\dots
\hat{x}_q\dots x_a}\.
\endsplit
$$
\endproclaim

\demo{Proof} Induction on $a$. For $a=1$ it is the Cartan formula. Suppose
that $a>1$. We have
$$
\align
di_{x_1\dots x_a}
&=d(i_{x_1}\cdot i_{x_2\dots x_a})=di_{x_1}\cdot
i_{x_2\dots x_a}-i_{x_1}\cdot di_{x_2\dots x_a}\\
\intertext{(by induction)}
&=x_{i_1}\cdot i_{x_2\dots x_a}-i_{x_1}\cdot \bigg(\sum_{p=2}^a(-1)^p
x_pi_{x_2\dots\hat{x}_p\dots x_a}\\
&\hskip100pt+\sum_{2\le p<q\le a}(-1)^{p+q}
i_{[x_p,x_q]x_2\dots\hat{x}_p\dots\hat{x}_q\dots x_a}\bigg)\\
\intertext{(we use $i_{x_1}x_p=x_pi_{x_1}+i_{[x_1,x_p]}$ and the
anticommutation of the various elements $i_x$)}
&=x_{i_1}\cdot i_{x_2\dots x_a}+\sum_{p=2}^a(-1)^{p-1}x_pi_{x_1\dots
\hat{x}_p\dots x_a}+\sum_{p=2}^a(-1)^{1+p}i_{[x_1,x_p]x_2\dots\hat{x}_p
\dots x_a}\\
&\qquad+\sum_{2\le p<q\le a}(-1)^{p+q}i_{[x_p,x_q]x_1\dots
\hat{x}_p\dots\hat{x}_q\dots x_a}\\
&=\sum_{p=1}^a(-1)^{p-1}x_pi_{x_1\dots\hat{x}_p\dots x_a}+
\sum_{1\le p<q\le a}(-1)^{p+q}i_{[x_p,x_q]x_1\dots\hat{x}_p\dots
\hat{x}_q\dots x_a}.
\endalign
$$
This proves the first equality. The second equality is proved in the same
manner, or may be deduced from the first one.\qed
\enddemo

\sbh{4.5.}Let $M$
be a $\frak{g}$-module. Recall that the complex
$C^{\bdot}(\frak{g};M)$ of cochains of $\frak{g}$ with
coefficients in $M$ is defined
by $C^a(\frak{g};M)=\op{Hom}\(\Lambda^a\frak{g},M)$;
the differential
$d\:\enspace C^{a-1}(\frak{g};M)\to C^{a}(\frak{g};M)$ acts as
$$
\split
d\phi(x_1\wedge\dots\wedge x_{a})
&=\sum_{p=1}^a(-1)^{p-1}x_p\phi
(x_1\wedge\dots\wedge\hat{x}_p\wedge\dots\wedge x_a)\\
&\qquad+
\sum_{1\le p<q\le a}(-1)^{p+q}\phi(x_1\wedge\dots\wedge\hat{x}_p\wedge
\dots\wedge\hat{x}_q\wedge\dots\wedge x_a)\.
\endsplit
$$
The differential in $C^{\bdot}(\frak{g};M)$ is called the
{\it Koszul differential}.

\remark{\rom{4.6.}\quad Remark} Consider $\U\frak{g}^{\bdot}$ as a
$\frak{g}$-module
by means of the left multiplication, where $\frak{g}$ is identified with
$\frak{g}^0$. Consider the
complex of cochains of $\frak{g}$
with coefficients in $\U\frak{g}^{\bdot}$,
$C^{\bdot}(\frak{g},\U\frak{g}^{\bdot})$. It is a double complex in which
the first differential is the Koszul differential,
and the second one is induced
by the differential in $\U\frak{g}^{\bdot}$. Let $C^{\bdot}$
be the associated total complex.
\endremark

For each $a\ge 0$,  we define an element $c^{a,-a}\in
C^a(\frak{g},\U\frak{g}^{-a})$ by
$$
c^{a,-a}(x_1\wedge\dots\wedge x_a)=i_{x_1\dots x_a},\qquad c^{00}=1\.
$$
One can reformulate the previous lemma as the following statement.

\itm{a}{The element $c=\sum_{a\ge 0} c^{a,-a}$
is a $0$-cocycle in $C^{\bdot}$.}

\sbh{4.7.}Suppose
we are given the collection of data (a)--(d) below.

\itm{a}{A Lie algebra $\frak{g}$;}

\itm{b}{a $\frak{g}$-module $M$;}

\itm{c}{a dg $\frak{g}^{\bdot}$-module $\Om^{\bdot}$;}

\itm{d}{an element $\om\in M\otimes\Om^n$, for some $n\ge 0$.}

\noindent Assume that the element
$\om$ satisfies the properties (i), (ii) below.

\itm{i}{$d\om=0$.}

\noindent Here $d=\op{Id}_M\otimes d_{\Om^{\bdot}}\:\enspace 
M\otimes\Om^n\to M\otimes\Om^{n+1}$.

\itm{ii}{$\om\in(M\otimes\Om^n)^{\frak{g}}$.}

\noindent Here the superscript $(\bcdot)^{\frak{g}}$
denotes the subspace of $\frak{g}$-invariants.
We consider $M\otimes\Om^n$ as a $\frak{g}$-module equal
to the tensor product
of the two $\frak{g}$-modules $M$ and $\Om^n$,
the last one being the $\frak{g}$-module
obtained via the identification $\frak{g}=\frak{g}^0$.

Consider the double complex
$C^{\bdot}(\frak{g};M\otimes\Om^{\bdot})$
defined by
$$
C^a(\frak{g};M\otimes\Om^b)=\op{Hom}\(\Lambda^a\frak{g},M\otimes\Om^b)\.
$$
The first differential is the Koszul differential in the standard
cochain complex of $\frak{g}$ with coefficients in $M\otimes\Om^{\bdot}$.
Here the action of $\frak{g}$ on $M\otimes\Om^{\bdot}$ is defined
{\it through the first factor\/} (i.e., as $x\cdot(m\otimes\a)=
xm\otimes\a$).
The second differential is induced by the differential in $\Om^{\bdot}$.
Let $C^{\bdot}(\frak{g};M,\Om^{\bdot})$ denote the associated
total complex.

For $0\le a\le n$,  we define the elements $w^a\in C^a(\frak{g};M\otimes
\Om^{n-a})$ by
$$
\om^0=0,\qquad
\om^a(x_1\wedge\dots\wedge x_a)=i_{x_1\dots x_a}\om\.
$$
Here the action of $\U\frak{g}^{\bdot}$ is defined through the second
factor.
Consider the element $\widehat{\om}=\sum_{a=0}^n\om^a\in
C^n(\frak{g};M,\Om^{\bdot})$.

\proclaim{\rom{4.8.}\quad Lemma}
The element $\widehat{\om}$ is an $n$-cocycle
in $C^{\bdot}(\frak{g};M,\Om^{\bdot})$.
\endproclaim

We shall call $\widehat{\om}$ the {\it Cartan cocycle\/} associated with
$\om$.

\demo{Proof} Let $d'$ (resp., $d''$) denote the first (resp., second)
differential in $C^{\bdot}(\frak{g};M\otimes\Om^{\bdot})$.
We have $d''(\om^0)=0$ by property 4.7(i). We have
$$
\split
&d'w^{a-1}(x_1\wedge\dots\wedge x_a)\\
&\qquad=\sum_{p=1}^a(-1)^{p-1}x_p\cdot
\om^{a-1}(x_1\wedge\cdots\wedge\hat{x}_p\wedge\cdots\wedge x_a)\\
&\qquad\qquad+\sum_{1\le p<q\le a}(-1)^{p+q}\om^{a-1}([x_p,x_q]\wedge x_1
\wedge\cdots\wedge\hat{x}_p\wedge\cdots\wedge\hat{x}_q\wedge\cdots
\wedge x_a)\,.
\endsplit
$$
Here $\omega\mapsto x^{(1)}\cdot\omega$,
(resp., $\omega\mapsto x^{(2)}\cdot\omega$),
 denotes the action of $\frak{g}$ on $M\otimes
\Om^n$ through the first (resp., second) factor.)

Using the identity
$$
\split
x_p\cdot(i_{x_1\dots\hat{x}_p\dots x_a}\om)&=
i_{x_1\dots\hat{x}_p\dots x_a}(x_p^{(1)}\om)\\
&=\text{(by 4.7(ii))}=-i_{x_1\dots\hat{x}_p\dots x_a}(x_p^{(2)}\om)=
-(i_{x_1\dots\hat{x}_p\dots x_a}x_p\)\om\.
\endsplit
$$

$$
\text{by Lemma 4.4 } = \,(di_{x_1\dots x_a}\)\om=
\text{(by (i))}=d(i_{x_1\dots x_a}\om)
=(d''\om^a)(x_1\wedge\cdots\wedge x_a).\qed
$$
\enddemo

\sbh{4.9.}Consider the dual complex $(\Om^{\bdot})^*$.
Note that the cocycle $\widehat{\om}$ defines a morphism of complexes
$$
\widehat{\om}\:\enspace (\Om^{n-\bdot})^*\to C^{\bdot}(\frak{g};M)\.
$$
Let us denote $H_i(\Om):=H^{-i}(\Om^*)$. The morphism $\widehat{\om}$
induces the following maps
$$
H^{-i}(\om)\:\enspace H_i(\Om^{\bdot})\to H^{n-i}(\frak{g};M)\qquad (0\le i\le n)\.
$$

\sbh{4.10.} The construction of this Section should be compared with 
~\cite{Br}. 

\head\S5.\enspace  Bosonization for the Virasoro algebra\endhead

Our account of the bosonization for the Virasoro algebra essentially
follows
the paper \cite{F}.

\sbh{5.1.}Heisenberg algebra

The {\it Heisenberg algebra\/} is the Lie algebra
$\Cal{H}$ defined by the generators $b_n$ (\<$n\in\Bbb Z$)
and $\bold{1}$, and relations

\itm{a}{$[b_n,b_m]=2n\delta_{n+m,0}\bold{1}$,
$[\bold{1},b_n]=0$ (\<$n,m\in\Bbb Z$).}

The elements $\bold{1}$ and $b_0$ lie in the center of $\Cal{H}$.
The elements
$b_n$, $n>0$, are called {\it annihilation operators}.

We shall denote by $\Cal{H}^+$
(resp., $\Cal{H}^-$, $\Cal{H}^0$) the Lie subalgebra
of $\Cal{H}$ generated by the elements $b_n$, $n>0$ (resp., by $b_n$,
$n<0$,
by $b_0$ and $\bold{1}$).

\sbh{5.2.}Given
two generators $b_n$, $b_m$, define their
{\it normal ordered product\/} $\Y{b_nb_m}\in\U\Cal{H}$ as follows.
If $n>0$ and $m\le 0$, we set $\Y{b_nb_m}=b_mb_n$; otherwise
set $\Y{b_nb_m} := b_nb_m$. We define 
$$
\{b_nb_m\}=b_nb_m-\Y{b_nb_m}
$$
(not to be confused with the Lie bracket!).

Similarly, the normal ordering of an arbitrary monomial
$\Y{b_{n_1}\cdots b_{n_a}}$ is defined: one
should pull all the annihilation
operators to the right, not changing their order. (The last requirement
does not matter very much since all the annihilation operators commute.)

We have by definition
$$
b_nb_m=\Y{b_nb_m}+\{b_nb_m\}\.
$$
The following identity generalizes this equality.

\proclaim{\rom{5.3.}\quad Wick theorem}
Let $c_1,\dots,c_a$, $c'_1,\dots,c'_b$ be
arbitrary elements of the set $\{b_n,\,n\in\Bbb Z\}$. We have
$$
\multline
\Y{c_1\cdots c_a}\,\cdot\,\Y{c'_1\cdots c'_b}
=\Y{c_1\cdots c'_b}+
\sum_{p,q}\{c_pc'_q\}\>\Y{c_1\cdots\hat{c}_p\cdots\hat{c}'_q\cdots c'_b}\\
+\sum_{p,p',q,q'}\{c_pc'_q\}\{c_{p'}c'_{q'}\}\>
\Y{c_1\cdots\hat{c}_p\cdots\hat{c}_{p'}\cdots\hat{c}'_q\cdots\hat{c}'_{q'}
\cdots c'_b}+\cdots.\qed
\endmultline
$$
\endproclaim

\sbh{5.4.}We introduce  a {\it ``free bosonic field''}, the formal expression:
$$
\phi(z)=-b_0\log\(z)+\sum_{n\ne 0}\frac{b_n}{n}\,z^{-n}\,,
$$
where $z$ is a formal variable.
Let $\phi'(z)$ be the formal derivative of $
\phi(z)$, that is the generating function
$$
\phi'(z)=-\sum_{n\in\Bbb Z} b_nz^{-n-1}.
$$
Here $z$ is a formal variable.

Commutation relations 5.1 (a) are equivalent to the identity
$$
\{\phi'(z\)\phi'(w)\}=\frac2{(z-w)^2}\.\tag"(a)"
$$

Let us deduce (a) from 5.1 (a). We have by definition
$$
\split
\{\phi'(z\)\phi'(w)\}
&=\phi'(z\)\phi'(w)-\Y{\phi'(z\)\phi'(w)}=
\sum_{n,m}(b_nb_m-\Y{b_nb_m}\)z^{-n-1}w^{-m-1}\\
&=\sum_{n>0}2nz^{-n-1}w^{n-1}=
2\paw \bigg(\sum_{n\ge 0}z^{-n-1}w^{n}\bigg)=
2z^{-1}\paw({1\over {(1-w/z)}})\\
&=2\paw \bigg(\frac1{z-w}\bigg)=\frac{2}{(z-w)^2}\.
\endsplit
$$

\sbh{5.5.}Given a complex number $\a\in\Bbb C$, define
the {\it Fock representation}
$F_\a $ of the Heisenberg algebra $\Cal{H}$ as an $\Cal{H}$-module
with one generator $v_\a $ and the relations
$$
b_nv_\a =0\quad (n>0),\qquad b_0v_\a =2\a v_\a ,\qquad
\bold{1} v_\a =v_\a \.
$$
The mapping $x\mapsto x\cdot v_\a $ gives an isomorphism of
$\Cal{H}^-$-modules $\U\Cal{H}^-@>\sim>>F_\a $.

Define the {\it shift operator}
$$
T_{\beta}\:F_\a \to F_{\a+\beta}
$$
as the unique $\Cal{H}^-$-linear operator sending $v_\a $ to
$v_{\a+\beta}$.

Define the {\it category $\Cal{O}$} to be
 the full subcategory of the category of
$\Cal{H}$-modules whose objects are representations $M$ having the
properties

\itm{a}{$\bold{1}$ acts as the identity on $M$;}

\itm{b}{$M$ is $b_0$-diagonalizable. All $b_0$-weight spaces
are finite-dimensional;}

\itm{c}{$M$ is $\Cal{H}^+$-locally finite. This
means that for any $x\in M$, the space $\U\Cal{H}^+\cdot x$ is
finite-dimensional.}

All Fock representations belong to the category $\Cal{O}$.

\sbh{5.6.}Recall that 
the {\it Witt algebra} $\Cal{L}$ is the Lie algebra of algebraic
vector fields on $\Bbb C^*=\Bbb C-\{0\}$. It can be defined by the
generators
$e_n=-z^{n+1}\,d/dz$ (\<$n\in\Bbb Z$) and the relations
$$
[e_n,e_m]=(n-m\)e_{n+m}\qquad (n,m\in\Bbb Z)\.
$$
The {\it Virasoro algebra} $\widehat{\Cal{L}}$ is
the Lie algebra defined by the generators
$L_n$ (\<$n\in\Bbb Z$) and $c$, and relations
$$
[L_n,L_m]=(n-m\)L_{n+m}+\frac{n^3-n}{12}\,\delta_{n+m,0}\cdot c,\quad
[c,L_n]=0\qquad (n,m\in\Bbb Z)\.
$$
We have the morphism of Lie algebras $\widehat{\Cal{L}}\to\Cal{L}$
sending $L_n$
to $e_n$, and $c$ to $0$. This morphism identifies $\Cal{L}$ with
$\widehat{\Cal{L}}/\Bbb C\cdot c$.
It identifies $\widehat{\Cal{L}}$ with a universal central extension
of $\Cal{L}$.

Let $\widehat{\Cal{L}}^+$ (resp., $\widehat{\Cal{L}}^-$,
$\widehat{\Cal{L}}^0$) denote the Lie subalgebras of $\widehat{\Cal{L}}$
generated by the elements $L_n$, $n>0$ (resp., by $L_n$,
$n<0$, by $L_0$, $c$).

The generating function
$$
T(z)=\sum_{n\in\Bbb Z} L_nz^{-n-2}
$$
is called the {\it stress-energy tensor}.

\sbh{5.7.}Let $\a_0\in\Bbb C$. Consider the expressions
$$
L_n=\frac14\sum_{p+q=n}\Y{b_pb_q}-\a_0(n+1\)b_n\qquad
(n\in\Bbb Z)\.\tag"(a)"
$$
Although they are infinite sums, these expressions are well defined
as operators on modules from the category $\Cal{O}$. We can rewrite them as
follows.
$$
\aligned
L_n&=\frac14\sum_{p\in\Bbb Z}b_{n-p}b_p-\a_0(n+1\)b_n\quad\text{if }\,
n\ne 0,\\
L_0&=\frac12\sum_{p\ge 1}b_{-p}b_p+\frac14b_0^2-\a_0b_0\.
\endaligned
\tag"(b)"
$$
In terms of generating functions, (a) reads
$$
T(z)=\frac14\,\Y{\phi'(z)^2}-\a_0\phi''(z)\.\tag"(c)"
$$

\proclaim{\rom{5.8.}\quad Theorem}
The expressions\/ \rom{5.7 (a)} define the
action of the Virasoro algebra on modules from the category $\Cal{O}$,
with the central charge $1-24\a_0^2$.
\endproclaim

\proclaim{\rom{5.9.}\quad Lemma} The operators $L_n$, \rom{5.7 (a)},
satisfy
the following commutation relations\/\rom:
$$
[b_n,L_m]=nb_{n+m}+2n(n-1\)\a_0\delta_{n,-m}\.\tag"(a)"
$$
\endproclaim

\demo{Proof} This can be checked easily using the definitions. Let us give
an alternative proof, using some simple {\it chiral calculus}.
We claim that (a) is equivalent to
$$
\phi'(z\)T(w)=-\frac{4\a_0}{(z-w)^3}
+\frac{\phi'(w)}{(z-w)^2}+\cdots\.\tag"(b)"
$$
Here (and below) dots denote the expression regular at $z=w$. Let us prove
that (b) implies (a). According to {\it Cauchy formula\/} of the chiral
calculus, we have
$$
[b_n,T(w)]=-\opl{res}_{z=w}(z^n\phi'(z\)T(w))
$$
for any $n$. By the binomial formula,
$$
\opl{res}_{z=w}\frac{z^n}{(z-w)^a}=\binom n{a-1}w^{n-a+1}.
$$
Therefore, (b) implies
$$
[b_n,T(w)]=2n(n-1\)\a_0w^{n-2}-n\phi'(w\)w^{n-1}\tag"(c)"
$$
which is equivalent to (a).

Let us prove (b). We have, by Wick's theorem,
$$
\gather
\phi'(z)\Y{\phi'(w)^2}=2\>\{\phi'(z\)\phi'(w)\}\>\phi'(w)+\cdots=
\frac{4\phi'(w)}{(z-w)^2}+\cdots,\\
\phi'(z\)\phi''(w)=\paw (\phi'(z\)\phi'(w))=\frac{4}{(z-w)^3}+\cdots\.
\endgather
$$
This implies (b) and proves the lemma.\qed
\enddemo

\subhead\rom{5.10.}\quad Proof of Theorem 5.8\endsubhead
The Virasoro commutation
relations can be verified directly, using the previous lemma and 5.7 (b),
by a tedious computation.

Let us give a proof using the chiral calculus.
We must prove the following.
$$
T(z\)T(w)=\frac{1-24\a_0^2}{2\(z-w)^4}+\frac{2T(w)}{(z-w)^2}+
\frac{T'(w)}{z-w}+\cdots\.\tag"(a)"
$$

Let us derive the Virasoro commutation relations from (a).
By the Cauchy formula of the chiral calculus, we have
$$
[L_n,T(w)]=\opl{res}_{z=w}(z^{n+1}T(z\)T(w))\.
$$
As in the proof of the previous lemma,
(a) implies that
$$
[L_n,T(w)]=T'(w\)w^{n+1}+2\(n+1\)T(w\)w^n+(1-24\a_0^2)\,
\frac{n^3-n}{12}\,w^{n-2}\tag"(b)"
$$
which is equivalent to the Virasoro commutation relations with the central
charge $1-24\a_0^2$.

Let us prove (a). We have, by Wick's theorem,
$$
\align
\Y{\phi'(z)^2}\Y{\phi'(w)^2}
&=2\>\{\phi'(z\)\phi'(w)\}^2+4\>\{\phi'(z\)\phi'(w)\}\>
\Y{\phi'(z\)\phi'(w)}\\
&=\frac{8}{(z-w)^4}+\frac{8}{(z-w)^2}\,\Y{\phi'(w)^2}+\frac{8}{z-w}
\Y{\phi'(w\)\phi''(w)}+\cdots,\\
\phi''(z)\Y{\phi'(w)^2}
&=2\>\{\phi''(z\)\phi'(w)\}\>\phi'(w)+\cdots\\
&=2\pa_z(\{\phi'(z\)\phi'(w)\}\)\phi'(w)=-\frac{8}{(z-w)^3}\,\phi'(w)+\cdots,\\
\Y{\phi'(z)^2}\phi''(w)
&=2\>\{\phi'(z\)\phi''(w)\}\>\phi'(z)=2\paw 
(\{\phi'(z\)\phi'(w)\}\)\phi'(z)\\
&=\frac{8}{(z-w)^3}\bigg(\phi'(w)+(z-w\)\phi''(w)+
\frac{(z-w)^2}{2}\,\phi'''(w)+\cdots\bigg)\\
&=\frac{8}{(z-w)^3}\,\phi'(w)+\frac{8}{(z-w)^2}\,\phi''(w)
+\frac{4}{z-w}\,\phi'''(w)+\cdots,\\
&\phi''(z\)\phi''(w)=\{\phi''(z\)\phi''(w)\}+\cdots
=\pa_z\paw (\{\phi'(z\)\phi'(w)\})+\cdots\\
&=-\frac{12}{(z-w)^4}+\cdots\.
\endalign
$$
Therefore we obtain, after adding,
$$
\align
T(z\)T(w)&=\bigg(\frac14\,\Y{\phi'(z)^2}-\a_0\phi''(z)\bigg)
\bigg(\frac14\,\Y{\phi'(w)^2}-\a_0\phi''(w)\bigg)\\
&=\frac{1-24\a_0^2}{2\(z-w)^4}
+\frac{(1/2)\Y{\phi'(w)^2}-2\a_0\phi''(w)}{(z-w)^2}\\
&\qquad+\frac{(1/2)\Y{\phi'(w\)\phi''(w)}-\a_0\phi'''(w)}{z-w}+\cdots\\
&=\frac{1-24\a_0^2}{2\(z-w)^4}+\frac{2T(w)}{(z-w)^2}+\frac{T'(w)}{z-w}
+\cdots\.
\endalign
$$
This proves (a) and the theorem. \qed

\sbh{5.11.}We define the representation $F_{\a;\a_0}$ of the
Virasoro algebra as the $\Cal{H}$-module $F_\a $, regarded
as an $\widehat{\Cal{L}}$-module by means of the
formulas in the previous theorem.

The representations $F_{\a,\a_0}$ will be called {\it Feigin--Fuchs
modules}.

\head\S6.\enspace  Bosonic vertex operators\endhead

\sbh{6.1.}
We need first to recall some formulas for 
"bosonization", see \cite{FF, TK}.

Let $\a,\beta\in\Bbb C$. Define the operators
$$
V_n(\beta)\:\enspace F_\a \to F_{\a+\beta}\qquad (n\in\Bbb Z)
$$
by means of the generating function
$$
V(\beta;z)=\sum_{n\in\Bbb Z}V_n(\beta\)z^{-n}\.
$$
We set by definition
$$
\split
V(\beta;z)&=T_{\beta}\Y{\!\exp\bigg(-\beta\sum_{n\ne 0}\frac{b_n}{n}\,
z^{-n}\bigg)}\\
&=T_{\beta}\exp\bigg(-\beta\sum_{n<0}\frac{b_n}{n}\,z^{-n}\bigg)
\exp\bigg(-\beta\sum_{n>0}\frac{b_n}{n}\,z^{-n}\bigg).
\endsplit
$$
This expression is an operator acting as follows
$$
V(\beta;z)\:\enspace F_\a \to F_{\a+\beta}((z^{-1}))\.
$$
These operators are called the {\it bosonic vertex operators}.
We also define the operators
$$
\align
V_-(\beta;z)&=\exp\bigg(-\beta\sum_{n<0}\frac{b_n}{n}\,z^{-n}\bigg)\:\enspace 
F_\a \to F_\a [[z]],\\
V_+(\beta;z)&=\exp\bigg(-\beta\sum_{n>0}\frac{b_n}{n}\,z^{-n}\bigg)\:\enspace 
F_\a \to F_\a [z^{-1}]\.
\endalign
$$
Thus, $V(\beta;z)=T_{\beta}V_-(\beta;z\)V_+(\beta;z)$.

\proclaim{\rom{6.2.}\quad Lemma}
We have
$$
\align
[b_n,V_-(\beta;z)]&=\cases2\beta z^nV_-(\beta;z)&\text{if $n>0$,}\\
0&\text{otherwise},\endcases\\
[b_n,V_+(\beta;z)]&=\cases2\beta z^nV_+(\beta;z)&\text{if $n<0$,}\\
0&\text{otherwise}.\endcases
\endalign
$$
\endproclaim

We leave the easy proof to the reader.

\proclaim{\rom{6.3.}\quad Theorem} For every $n\in\Bbb Z$, we have
$[b_n,V(\beta;z)]=2\beta z^nV(\beta;z)$.
\endproclaim

\demo{Proof} Follows at once from the previous lemma.\qed
\enddemo

\sbh{6.4.}Let us give an alternative proof, which uses the chiral calculus.
Let us introduce the expression
$$
\wt{V}(\beta;z)=z^{2\beta\a}V(\beta;z)\.
$$
At this point $z^{\beta\a}$ is a formal symbol. It will play its role
in the sequel, when we start differentiating. More precisely,
$\wt{V}(\beta;z)$ is the operator $V(\beta;z)$ considered as a section
of a certain De Rham complex with a nontrivial connection.

Our formula is equivalent to
$[b_n,\wt{V}(\beta;z)]=2\beta\wt{V}(\beta;z)$.
We must prove that
$$
\phi'(z\)\wt{V}(\beta;w)=-\frac{2\beta}{z-w}\,
\wt{V}(\beta;w)+\cdots\.\tag"(a)"
$$

Recall  the {\it ``free bosonic field''}
$$
\phi(z)=q-b_0\log\(z)+\sum_{n\ne 0}\frac{b_n}{n}\,z^{-n}\.
$$
where now we added a $z$-independent term, an operator $q$
that satisfies the following commutation relations
$$
[q,b_n]=2\delta_{n,0}\bold{1}\.
$$
It follows that
$$
[b_0,e^{-\beta q}]=2\beta e^{-\beta q}.
$$
We identify the operator $e^{-\beta q}$ with $T_{\beta}$. Thus,
$$
\wt{V}(\beta;z)=\Y{\!\exp(-\beta\phi(z))},
$$
where the only nontrivial normal ordering with the operator $q$ is defined
as
$\Y{b_0q}=qb_0$.

We have
$$
\phi(z\)\phi(w)=2\log\(z-w)+\Y{\phi(z\)\phi(w)}\.\tag"(b)"
$$
It follows that
$$
\phi'(z\)\phi(w)=\frac{2}{z-w}+\cdots,
$$
hence, from the Wick theorem
$$
\phi'(z)\Y{\phi(w)^n}=\frac{2n}{z-w}\,\Y{\phi(w)^{n-1}}+\cdots
$$
(\<$n\ge 0$), so if $f(\phi)$ is any power series in $\phi$,
$$
\phi'(z)\Y{f(\phi(w))}=\frac{2}{z-w}\,\Y{f'(\phi(w))}+\cdots\.
$$
In particular, we have
$$
\phi'(z)\Y{\!\exp\(-\beta\phi(w))}=-\frac{2\beta}{z-w}
\Y{\!\exp\(-\beta\phi(w))}+\cdots,
$$
which proves (a) and the theorem.\qed

\sbh{6.5.}We  regard the operators $V_n(\beta)$ as operators acting
on Feigin--Fuchs modules
$$
V_n(\beta)\:\enspace F_{\a;\a_0}\to F_{\a+\beta;\a_0}\.
$$
The generating function
$V(\beta;z)$ will be understood in the same sense.

\proclaim{\rom{6.6.}\quad Theorem}
For any $n\in\Bbb Z$, we have
$$
[L_n,V(\beta;z)]=\bigg(z^{n+1}\,\frac{d}{dz}+((\beta^2-2\a_0\beta)
(n+1)+2\a\beta\)z^n\bigg)V(\beta;z)\.
$$
In other words,
$$
[L_n,\wt{V}(\beta;z)]
=\bigg(z^{n+1}\,\frac{d}{dz}+(\beta^2-2\a_0\beta)(n+1\)z^n\bigg)
\wt{V}(\beta;z)\.
$$
\endproclaim

\demo{Proof} We must prove that
$$
T(z\)\wt{V}(\beta;w)=\frac{\beta^2-2\a_0\beta}{(z-w)^2}\,
\wt{V}(\beta;w)+\frac{1}{z-w}\,\wt{V}'(\beta;w)+\cdots\.
\tag"(a)"
$$
Let us prove (a). We have, by the Wick theorem,
$$
\Y{\phi'(z)^2}\Y{\phi(w)^n}=\frac{4n(n-1)}{(z-w)^2}\Y{\phi(w)^{n-2}}+
\frac{4n}{z-w}\,\Y{\phi'(w)\phi(w)^{n-1}}+\cdots
$$
for any $n\ge 0$, hence
$$
\Y{\phi'(z)^2}\wt{V}(\beta;w)
=\frac{4\beta^2}{(z-w)^2}\,\wt{V}(\beta;w)+\frac{4}{z-w}\,
\wt{V}'(\beta;w)+\cdots.
$$
It follows from 6.4 (a) that
$$
\phi''(z\)\wt{V}(\beta;w)
=\frac{2\beta}{(z-w)^2}\,\wt{V}(\beta;w)+\cdots\.
$$
Summing, we get (a). This proves the theorem.\qed
\enddemo

\proclaim{\rom{6.7.}\quad Lemma} We have
$$
V_+(\beta_1;z_1\)V_-(\beta_2;z_2)=(1-z_1^{-1}z_2)^{2\beta_1\beta_2}
V_-(\beta_2;z_2\)V_+(\beta_1;z_1)
$$
\rom(equality in $\op{Hom}\(F_\a ,F_\a [[z_1^{-1},z_2]])$\<\rom).
\endproclaim

Here we understand $(1-z_1^{-1}z_2)^{2\beta_1\beta_2}$ as the formal power
series
$$
\exp\bigg(-2\beta_1\beta_2\sum_{n>0}z_1^{-n}\frac{z_2^n}n\bigg).
$$

\remark{Remark} The operator $V_-(\beta_2;z_2\)V_+(\beta_1;z_1)$
belongs to the space $\op{Hom}\(F_\a ,\,F_\a [z_1^{-1},$
$z_2])$.
\endremark

\demo{Proof} By Lemma 6.2, we have
$$
V_+(\beta_1;z_1\)b_{-n}=(b_{-n}-2\beta_1z_1^{-n}\)V_+(\beta_1;z_1)
$$
for $n>0$, hence
$$
\split
V_+(\beta_1;z_1)\exp\bigg(\beta_2\,\frac{b_{-n}}{n}\,z_2^n\bigg)
&=\exp\bigg(\beta_2\,\frac{b_{-n}-2\beta_1z_1^{-n}}{n}\,z_2^n\bigg)
V^+(\beta_1;z_1)\\
&=\exp\bigg(-2\beta_1\beta_2\,\frac1n\,z_1^{-n}z_2^n\bigg)
\exp\bigg(\beta_2\,\frac{b_{-n}}{n}\,z_2^n\bigg)V_+(\beta_1;z_1)
\endsplit
$$
Therefore,
$$
\split
V_+(\beta_1;z_1\)V_-(\beta_2;z_2)
&=\exp\bigg(-2\beta_1\beta_2\sum_{n>0}
\frac{1}{n}z_1^{-n}z_2^n\bigg)V_-(\beta_2;z_2\)V_+(\beta_1;z_1)\\
&=\exp\(2\beta_1\beta_2\log\(1-z_1^{-1}z_2)\)V_-(\beta_2;z_2\)
V_+(\beta_1;z_1)\\
&=(1-z_1^{-1}z_2)^{2\beta_1\beta_2}V_-(\beta_2;z_2\)V_+(\beta_1;z_1)\.\qed
\endsplit
$$
\enddemo

\sbh{6.8.}Choose complex numbers $\a,\beta_1,\dots,\beta_p$.
Consider the generating function for all  compositions
of the form
$$
F_\a@>V_{n_p}(\beta_p)>>F_{\a+\beta_p}@>>>
\cdots@>V_{n_1}(\beta_1)>>F_{\a+\beta_p+\dots+\beta_1};
$$
We have an equation
$$
\sum_{(n_1,\dots,n_p)\in\Bbb Z^p}V_{n_1}(\beta_1)\cdots
V_{n_p}(\beta_p\)z_1^{-n_1}\cdots z_p^{-n_p}=
V(\beta_1;z_1)\cdots V(\beta_p;z_p)
$$
as a formal power series. The previous lemma shows that
$$
V(\beta_1;z_1)\cdots V(\beta_p;z_p)=\prod_{1\le i<j\le p}
(1-z_i^{-1}z_j)^{2\beta_i\beta_j}\Y{V(\beta_1;z_1)\cdots V(\beta_p;z_p)}\.
\tag"(a)"
$$

\sbh{6.9.}Let us look more closely at the operator
$\Y{V(\beta_1;z_1)\cdots V(\beta_p;z_p)}$.

We have
$$
\split
\Y{V(\beta_1;z_1)\cdots V(\beta_p;z_p)}
&=T_{\beta_1+\cdots+\beta_p}\exp\bigg(-\sum_{n<0}
\frac{\beta_1z_1^{-n}+\cdots+\beta_pz_p^{-n}}{n}\,b_n\bigg)\\
&\qquad\times
\exp\bigg(-\sum_{n>0}\frac{\beta_1z_1^{-n}+\cdots+\beta_pz_p^{-n}}{n}\,
b_n\bigg).
\endsplit
$$
It follows that

\itm{a}{the operator $\Y{V(\beta_1;z_1)\cdots V(\beta_p;z_p)}$
belongs to the space
$\op{Hom}\(F_\a ,\,F_{\a+\sum\beta_i}[z_1,$
$z_1^{-1},\dots,z_p,z_p^{-1}])$.}

\noindent In particular, this operator is a holomorphic
operator-valued function on
the space $\Bbb C^p-\bigcup_{i=1}^p\{z_i=0\}$.

Also, the above formula shows that

\itm{b}{for any bijection $\sigma\:\{1,\dots,p\}@>\sim>>\{1,\dots,p\}$,
we have}
$$
\Y{V(\beta_1;z_1)\cdots V(\beta_p;z_p)}=
\Y{V(\beta_{\sigma(1)};z_{\sigma(1)})\cdots V(\beta_{\sigma(p)}
z_{\sigma(p)})}\.
$$

\sbh{6.10.}Formula 6.8 (a) shows that the formal power series
$V(\beta_1;z_1)\cdots V(\beta_p;z_p)$ defines the germ of a holomorphic
multivalued function in the domain $|z_1|>\cdots>|z_p|>0$,
where $(1-z_i^{-1}z_j)^{2\beta_i\beta_j}$ is understood as
$\exp\(-2\beta_i\beta_j\sum_{n>0}z_i^{-n}z_j^n/n)$.

\sbh{6.11.}Let us introduce the tilded operators. We set
$$
\wt{V}(\beta_i;z_i)=\Y{\!\exp\(-\beta_i\phi(z_i))}=
T_{\beta_i}z_i^{\beta_ib_0}\Y{\!\exp\bigg(-\beta_i\sum_{n\ne 0}
\frac{b_n}{n}\,z_i^{-n}\bigg)}\.
$$
Since
$z_i^{\beta_ib_0}T_{\beta_j}=z^{2\beta_i\beta_j}T_{\beta_j}z_i^{\beta_i
b_0}$, we have
$$
\split
\wt{V}(\beta_1;z_1)\cdots\wt{V}(\beta_p;z_p)
&=T_{\sum\beta_i}\prod_iz_i^{\beta_ib_0}\prod_{i<j}z_i^{2\beta_j\beta_j}
V(\beta_1;z_1)\cdots V(\beta_p;z_p)\\
&=\prod_i z_i^{2\beta_i\a}\prod_{i<j}(z_i-z_j)^{2\beta_i\beta_j}
\Y{V(\beta_1;z_1)\cdots V(\beta_p;z_p)}\.
\endsplit
$$
Recall that the operator $\Y{V(\beta_1;z_1)\cdots V(\beta_p;z_p)}$
belongs to the space
$$
\op{Hom}\(F_\a ,F_{\a+\sum\beta_i}[z_1,z_1^{-1},\dots,z_p,z_p^{-1}])\.
$$
This defines
the product $\wt{V}(\beta_1;z_1)\cdots\wt{V}(\beta_p;z_p)$ as the germ
of the multivalued holomorphic function in the domain
$$
\{(z_1,\dots,z_p)\in\Bbb C^p\mid z_i\notin\Bbb R_{\le 0}\,
\text{ for all }\,i;\,|z_1|>\cdots>|z_p|\}\.
$$
Here $z^{\ga}$ is understood as $\exp\(\ga\log\(z))$,
where $\log\(z)$
is the branch of the logarithm that takes real values for $z\in\Bbb R_{>0}$.

Let us regard the previous compositions as operators acting on the
Feigin--Fuchs modules
$$
\wt{V}(\beta_1;z_1)\cdots\wt{V}_p(\beta_p;z_p)\:\enspace 
F_{\a;\a_0}\to F_{\a+\sum\beta_i;\a_0}\.
$$

\proclaim{\rom{6.12.}\quad Theorem}
For any $n\in\Bbb Z$,
$$
\split
&[L_n,\wt{V}(\beta_1;z_1)\cdots\wt{V}(\beta_p;z_p)]\\
&\qquad=
\bigg(\sum_{i=1}^pz_i^{n+1}\pa_{z_i}+(\beta_i^2-2\a_0\beta_i)(n+1\)z_i^n
\bigg)\wt{V}(\beta_1;z_1)\cdots\wt{V}(\beta_p;z_p)\.
\endsplit
$$
In other words,
$$
\split
&[L_n,\Y{V(\beta_1,z_1)\cdots V(\beta_p;z_p)}]\\
&\qquad=\botsmash{\bigg(\sum_{i=1}^p}
(z_i^{n+1}\pa_{z_i}+((\beta_i^2-2\a_0\beta_i)(n+1)+2\beta_i\a\)z_i^n)\\
&\hskip110pt+\sum_{1\le i<j\le p}2\beta_i\beta_j\frac{z_i^{n+1}-z_j^{n+1}}
{z_i-z_j}\bigg)\Y{V(\beta_1;z_1)\cdots V(\beta_p;z_p)}\.
\endsplit
$$
\endproclaim

\demo{Proof} By the Cauchy formula,
$$
[L_n,\wt{V}(\beta_1;z_1)\cdots\wt{V}(\beta_p;z_p)]=
\sum_{i=1}^p\opl{res}_{z=z_i}(z^{n+1}T(z\)\wt{V}
(\beta_1;z_1)\cdots\wt{V}_p(\beta_p;z_p))\.
$$
Note that
$$
\wt{V}(\beta_1;z_1)\cdots\wt{V}(\beta_p;z_p)
=\prod_{i<j}(z_i-z_j)^{2\beta_i\beta_j}\Y{\wt{V}(\beta_1;z_1)
\cdots\wt{V}(\beta_p;z_p)}\.
$$
We claim that
$$
\split
&T(z)\;\Y{\wt{V}(\beta_1;z_1)\cdots\wt{V}(\beta_p;z_p)}\\
&\quad=\bigg(\sum_{i=1}^p\frac{\beta_i^2-2\a_0\beta_i}{(z-z_i)^2}+
\sum_{i<j}\frac{2\beta_i\beta_j}{(z-z_i)(z-z_j)}\bigg)
\;\Y{\wt{V}(\beta_1;z_1)\cdots\wt{V}(\beta_p;z_p)}\\
&\qquad\qquad+\sum_{i=1}^p\frac{1}{z-z_i}\;
\Y{\wt{V}(\beta_1;z_1)\cdots\wt{V}'(\beta_i;z_i)
\cdots\wt{V}(\beta_p;z_p)}+\cdots\.
\endsplit
\tag"(a)"
$$

Let us prove (a). To shorten the formulas, assume that $p=2$; the general
case
is completely similar. By Wick's theorem, we have
$$
\split
&\,\Y{\phi'(z)^2}\,\Y{\phi(z_1)^{n_1}\phi(z_2)^{n_2}}\\
&\quad=\frac{4n_1(n_1-1)}{(z-z_1)^2}\,
\,\Y{\phi(z_1)^{n_1-2}\phi(z_2)^{n_2}}
+\frac{4n_2(n_2-1)}{(z-z_2)^2}\,
\,\Y{\phi(z_1)^{n_1}\phi(z_2)^{n_2-2}}\\
&\qquad+\frac{8n_1n_2}{(z-z_1)(z-z_2)}\,
\,\Y{\phi(z_1)^{n_1-1}\phi(z_2)^{n_2-1}}+
\frac{4n_1}{z-z_1}\,\,\Y{\phi'(z_1\)\phi(z_1)^{n_1-1}\phi(z_2)^{n_2}}\\
&\qquad+\frac{4n_2}{z-z_2}\,
\,\Y{\phi(z_1)^{n_1}\phi'(z_2\)\phi(z_2)^{n_2-1}}+\cdots
\endsplit
$$
for any $n_1,n_2\ge 0$, hence
$$
\split
\,\Y{\phi'(z)^2}\,\Y{\wt{V}(\beta_1;z_1)\wt{V}(\beta_2;z_2)}
&=4\bigg(\frac{\beta_1}{z-z_1}+\frac{\beta_2}{z-z_2}\bigg)^2
\,\Y{\wt{V}(\beta_1;z_1\)\wt{V}(\beta_2;z_2)}\\
&\qquad+4\bigg(\frac{\pa_{z_1}}{z-z_1}+\frac{\pa_{z_2}}{z-z_2}\bigg)
\,\Y{\wt{V}(\beta_1;z_1\)\wt{V}(\beta_2;z_2)}+\cdots\.
\endsplit
$$
Similarly, the application of the Wick theorem gives
$$
\phi'(z)\,\Y{\wt{V}(\beta_1;z_1\)\wt{V}(\beta_2;z_2)}=
\bigg(-\frac{2\beta_1}{z-z_1}-\frac{2\beta_2}{z-z_2}\bigg)
\,\Y{\wt{V}(\beta_1;z_1\)\wt{V}(\beta_2;z_2)}+\cdots,
$$
whence
$$
\phi''(z)\,\Y{\wt{V}(\beta_1;z_1\)\wt{V}(\beta_2;z_2)}
=\bigg(\frac{2\beta_1}{(z-z_1)^2}+
\frac{2\beta_2}{(z-z_2)^2}\bigg)\,\Y{\wt{V}(\beta_1;z_1\)\wt{V}(\beta_2;z_2)}
+\cdots\.
$$
Summing, we get formula (a).

The theorem follows from formula (a) by the above mentioned
Cauchy residue formula. One should take into account that
$$
\Big(\opl{res}_{z=z_i}+\opl{res}_{z=z_j}\Big)\frac{z^{n+1}}
{(z-z_i)(z-z_j)}=
\frac{z_i^{n+1}-z_j^{n+1}}{z_i-z_j}
$$
and
$$
(z_i^{n+1}\pa_{z_i}+z_j^{n+1}\pa_{z_j})(z_i-z_j)^{2\beta_i\beta_j}=
2\beta_i\beta_j\,\frac{z_i^{n+1}-z_j^{n+1}}{z_i-z_j}\,
(z_i-z_j)^{2\beta_i\beta_j}.
$$
This completes the proof.\qed
\enddemo

\head\S7.\enspace  Screening charges\endhead

\sbh{7.1.}We fix the parameters $\a,\a_0\in\Bbb C$ as in the
previous section. Let $\beta_+$, $\beta_-$ be the complex numbers defined
by
$$
\beta_{\pm}=\a_0\pm\sqrt{\a_0+1}\.
$$
Thus, $\beta_{\pm}$ are the two roots of the equation
$\beta^2-2\a_0\beta=1$.

The vertex operators $V(\beta_+;z)$, $V(\beta_-;z)$ are called the
{\it
screening charges}. By theorem 6.6, for any $n\in\Bbb Z$,
$$
[L_n,V(\beta_{\pm};z)]=\bigg(\frac{d}{dz}+\frac{2\a\beta_{\pm}}{z}\bigg)
(z^{n+1}V(\beta_{\pm};z))\.\tag"(a)"
$$
In other words,
$$
T(z\)\wt{V}(\beta_{\pm};z)=\frac{\wt{V}(\beta_{\pm};w)}{(z-w)^2}+
\frac{\wt{V}'(\beta_{\pm};w)}{z-w}+
\cdots=\paw \bigg(\frac{\wt{V}(\beta_{\pm};w)}{z-w}\bigg)+\cdots.\tag"(b)"
$$

\sbh{7.2.}Let us introduce the operators
$$
T(L_n;z)=z^{n+1}T(z),\quad V_{\pm}(e_n;z)=z^n\wt{V}(\beta_{\pm};z)\qquad
(L_n\in\widehat{\Cal{L}},\;e_n\in\Cal{L})\.
$$
It follows from (b) above that these operators satisfy the following
cocycle property
$$
T(L_n;z\)V_{\pm}(e_m;w)-T(L_m;z\)V_{\pm}(e_n;w)=
\frac{V_{\pm}([e_n,e_m];w)}{z-w}+\cdots\.\tag"(a)"
$$

\sbh{7.3.}More generally, suppose that $\beta_1,\dots,\beta_p$ is a
sequence
of complex numbers, each $\beta_i$ being equal to $\beta_-$ or $\beta_+$.
Let us call such a sequence a {\it screening sequence}.

Let us
consider the operator
$\,\Y{V(\beta_1;z_1)\cdots V(\beta_p;z_p)}$. Note that by the symmetry
property 6.9 (b), we may actually assume that
$\beta_1=\cdots=\beta_{p'}=\beta_-$ and $\beta_{p'+1}=\cdots=\beta_p=
\beta_+$.

It follows from Theorem 6.12, that
$$
\split
&[L_n,\,\Y{V(\beta_1;z_1)\cdots V(\beta_p;z_p)}]\\
&\qquad=\sum_{i=1}^p\bigg(\pa_{z_i}+\frac{2\a\beta_i}{z_i}+
\sum_{j\ne i}\frac{2\beta_i\beta_j}{z_i-z_j}\bigg)(z_i^{n+1}
\,\Y{V(\beta_1;z_1)\cdots V(\beta_p;z_p)})
\endsplit
$$
for any screening sequence $\beta_1,\dots,\beta_p$, and
for all $n\in\Bbb Z$.

In the sequel we fix the screening sequence
$\beta_1=\cdots=\beta_p=\beta_+$.
However, all the constructions below are valid, with the obvious
modifications, for an arbitrary screening sequence.

\sbh{7.4.}Consider the ring
$$
A_p=\Bbb C[[z_1,\dots,z_p]]
\bigg[\prod_{i=1}^pz_i^{-1},\prod_{1\le i<j\le p}(z_i-z_j)^{-1}\bigg]
$$
(thus $A_p$ is the ring of functions on the formal variety $X_p$, the
$p$\<th power of the formal punctured disk without diagonals).

For $1\le a\le p$, let $\Om^a$ denote the space of algebraic
differential $a$-forms on $X_p$. Thus, $\Om^0=A_p$, and the elements of
$\Om^a$ have
the form
$$
\sum f(z_1,\dots,z_p)\,dz_{i_1}\wedge\dots\wedge dz_{i_a}\qquad
(f(z_1,\dots,z_p)\in A_p)\.
$$

Consider the complex
$$
\Om^{\bdot}\:\enspace 0@>>>\Om^0@>d>>\Om^1@>d>>
\cdots@>d>>\Om^p@>>>0,
$$
where the differential is defined by
$$
d\eta=d_{DR}(\eta)+\bigg(\sum_{i=1}^p2\a\beta_+\,\frac{dz_i}{z_i}+
\sum_{1\le i<j\le p}2\beta_+^2\,\frac{dz_i-dz_j}{z_i-z_j}\bigg)\wedge\eta
$$
(\<$\eta\in\Om^a$), $d_{DR}$ being the usual de Rham differential.

Each vector field $\tau=\mu(z\)\pa_z\in\Cal{L}$,
$\mu(z)\in\Bbb C[z,z^{-1}]$,
defines a series of morphisms $i_{\tau}\:\Om^a\to\Om^{a-1}$, by
$$
\split
&i_{\tau}(f(z_1,\dots,z_p)\,dz_{i_1}\wedge\dots\wedge dz_{i_a})\\
&\qquad=
\sum_{b=1}^a(-1)^{b-1}\mu(z_{i_b}\)f(z_1,\dots,z_p)\,dz_{i_1}\wedge\cdots
\wedge\hat{dz}_{i_b}\wedge\cdots\wedge dz_{i_a}\.
\endsplit
$$
Let us define the morphisms of
the {\it Lie derivative} $\op{Lie}_{\tau}\:\Om^a\to\Om^a$ by
$$
\op{Lie}_{\tau}(\eta)=di_{\tau}(\eta)+i_{\tau}(d\eta)\.
$$
The morphisms $i_{\tau},\op{Lie}_{\tau}$ define the action of the dg
Lie algebra $\Cal{L}^{\bdot}$ associated with the Witt algebra $\Cal{L}$
(cf.~4.2),
on the complex $\Om^{\bdot}$ (cf. 4.2).

\sbh{7.5.}Consider the complex $\op{Hom}\(F_{\a;\a_0},
F_{\a+p\beta_+;\a_0}\otimes\Om^{\bdot})$.

Let us note that
the action of the Virasoro algebra $\widehat{\Cal{L}}$ on the modules
$F_{\a;\a_0}$, $F_{\a+p\beta_+;\a}$ induces the action of
$\widehat{\Cal{L}}$ on the space
$\op{Hom}\(F_{\a;\a_0},F_{\a+p\beta^+;\a_0})$
by the usual commutator formula, which factors through $\Cal{L}$, since
the Feigin--Fuchs modules have the same central charge. This in turn
induces the action of $\Cal{L}$ on the complex
$\op{Hom}\(F_{\a;\a_0},F_{\a+p\beta_+;\a_0}\otimes\Om^{\bdot})$,
which we shall call the {\it first action}.
We have also the {\it second action} of $\Cal{L}$ on
$\op{Hom}\(F_{\a;\a_0},F_{\a+p\beta_+;\a_0}\otimes
\Om^{\bdot})$ which is induced by the action of $\Cal{L}$ on
$\Om^{\bdot}$ through the Lie derivative. The first and the second
actions commute.

We also have the operators
$$
i_{\tau}\:\enspace \op{Hom}\(F_{\a;\a_0},F_{\a+p\beta_+;\a_0}\otimes\Om^{\bdot})
\to\op{Hom}\(F_{\a;\a_0},F_{\a+p\beta_+;\a_0}\otimes\Om^{\bdot})[-1]\qquad
(\tau\in\Cal{L})
$$
induced by the operators of the same name
on $\Om^{\bdot}$. These operators,
together with the second action, define the action of the dg Lie algebra
$\Cal{L}^{\bdot}$ on the complex
$\op{Hom}\(F_{\a;\a_0},F_{\a+p\beta_+;\a_0}\otimes\Om^{\bdot})$.

Let us define the element
$V^{0p}=V^{0p}(z_1,\dots,z_p)\in\op{Hom}\(F_{\a;\a_0},
F_{\a+p\beta_+;\a_0}\otimes\Om^{p})$ by
$$
V^{0p}(z_1,\dots,z_p)=\,\Y{V(\beta_+;z_1)\cdots V(\beta_+;z_p)}
\,dz_1\wedge\cdots\wedge dz_p\.
$$
The formula in~7.3 may be reformulated as the following theorem.

\proclaim{\rom{7.6.}\quad Theorem}
The element $V^{0p}$ lies in the subspace of $\Cal{L}$-invariants
$$
\op{Hom}\(F_{\a;\a_0},F_{\a+p\beta_+;\a_0}\otimes\Om^{p})^{\Cal{L}}.
$$
Here the action of the Lie algebra $\Cal{L}$ on the space
$\op{Hom}\(F_{\a;\a_0},F_{\a+p\beta_+;\a_0}\otimes\Om^{p})$
is the sum of the first and the second actions.\qed
\endproclaim

\sbh{7.7.}Now  we apply the construction of \S4. Consider
the double complex
$$
C^{\bdot}(\Cal{L};\op{Hom}\(F_{\a;\a_0},F_{\a+p\beta_+;\a_0}\otimes
\Om^{\bdot}))\.
$$
Here the first differential is the Koszul differential of the cochain
complex
of the Lie algebra $\Cal{L}$ with coefficients in the module
$\op{Hom}\(F_{\a;\a_0},F_{\a+p\beta_+;\a_0}\otimes
\Om^{\bdot})$,
with the first $\Cal{L}$-action. The second differential is induced by the
differential in $\Om^{\bdot}$.

Define the elements
$$
V^{a,p-a}\in C^a(\Cal{L};\op{Hom}\(F_{\a;\a_0},F_{\a+p\beta_+;\a_0}
\otimes\Om^{p-a}))\qquad (0\le a\le p)
$$
by
$$
V^{a,p-a}(\tau_1\wedge\cdots\wedge\tau_a)=i_{\tau_1}\cdots i_{\tau_a}
(V^{0p})\.
$$

\proclaim{\rom{7.8.}\quad Theorem}
The element $V=(V^{0p},\dots,V^{p0})$ is
a $p$-cocycle in the total complex associated with the double complex
$C^{\bdot}(\Cal{L};\op{Hom}\(F_{\a;\a_0},F_{\a+p\beta_+;\a_0}\otimes
\Om^{\bdot}))$.
\endproclaim

The proof is the same as that of Lemma 4.8.

\proclaim{\rom{7.9.}\quad Corollary} The cocycle $V$ induces the maps
$$
f_i\:\enspace H^i(\Om^{\bdot})^*\to\op{Ext}^{p-i}_{\widehat{\Cal{L}}}(F_{\a;\a_0},
F_{\a+p\beta_+;\a_0}),\qquad i=0,\dots,p\.\qed
$$
\endproclaim

The map
$$
f_p\:\enspace H^p(\Cal{L})^*\to\op{Hom}_{\widehat{\Cal{L}}}(F_{\a;\a_0},F_{\a+p\beta_
+;
\a_0})
$$
is the Feigin--Fuchs intertwiner, \cite{FF, Chapter 4}.

\head Chapter 3. $\widehat{\frak{sl}}(2)$-case\endhead

The main construction of this and the next chapters was inspired by
\cite{BMP1, BMP2}.

\head\S8.\enspace  Wakimoto realization\endhead

\sbh{8.1.}Let $\frak{g}$ be a Lie algebra and $B(\bcdot,\bcdot)$
an invariant bilinear form on $\frak{g}$.
The corresponding {\it affine Lie algebra}
$\hat{\frak{g}}$ is defined by the generators
$X_n$ (\<$X\in\frak{g}$, $n\in\Bbb Z$\<) and $\bold{1}$,
and the relations
$$
[X_n,X_m]=[X,Y]_{m+n}+nB(X,Y\)\delta_{m+n,0}\cdot\bold{1}\qquad
(X,Y\in\frak{g},\;m,n\in\Bbb Z)\.\tag"(a)"
$$
The element $\bold{1}$ will act as the identity on all our representations.

Let us introduce the generating functions ({\it currents\/})
$X(z)=\sum_{n\in\Bbb Z} X_nz^{-n-1}$ (\<$X\in\frak{g}$).
Formula (a) is equivalent
to
$$
X(z\)Y(w)=\frac{B(X,Y)}{(z-w)^2}+\frac{[X,Y](w)}{z-w}+\cdots\.\tag"(b)"
$$
Here as usual the dots stands for the
part regular at $z=w$. One deduces (a) from
(b) at once, using the chiral Cauchy formula
$$
[X_n,Y(w)]=\opl{res}_{z=w}(z^nX(z\)Y(w))\.
$$

\sbh{8.2.}In this chapter we assume that $\frak{g}=\frak{sl}(2)$
with the standard generators $E$, $F$, $H$.
For $X,Y\in\frak{g}$, we set $(X,Y)=\op{tr}\(XY)$.
Thus, $(E,F)=(F,E)=1$, $(H,H)=2$.
We fix a complex number $k$, and set $B(X,Y)=k(X,Y)$.

The bosonization formulas for the algebra $\widehat{\frak{g}}$
presented below were discovered by M.~Wakimoto, \cite{W}.

\sbh{8.3.}Let $\bold{a}$ denote the Lie algebra defined by the generators
$b_n$, $a_n$, $a^*_n$ (\<$n\in\Bbb Z$),
$\bold{1}$ and the relations

\itm{a}{$[b_n,b_m]=2n\delta_{n+m,0}\cdot\bold{1}$;}

\itm{b}{$[a_n,a^*_m]=\delta_{n+m,0}\cdot\bold{1}$,
$[a_n,a_m]=[a^*_n,a^*_m]=0$;}

\itm{c}{$[b_n,a_m]=[b_n,a^*_m]=0$, $\bold{1}$ commutes with everything.}

Let $\bold{a}_+$ denote the Lie subalgebra of $\bold{a}$ generated by the
elements
$b_n$, $a^*_n$ (\<$n>0$), $a_n$ (\<$n\ge 0$).
These generators are called
{\it annihilation operators}. One introduces the {\it normal ordering}
of a monomial in $\U\bold{a}$
in the usual way: all annihilation operators should be pulled to the right.

All $\bold{a}$-modules $M$ that we
consider be $\bold{a}_+$-{\it locally finite}, i.e.,
will have the property
that for every $x\in M$, the space $(\U\bold{a}_+)x$ is finite-dimensional.
The operator $\bold{1}$ will act as the identity.

For computational purposes, we shall also use one more operator $q$,
with the only nontrivial commutation relation
$$
[q,b_0]=\bold{1}\.
$$

\sbh{8.4.}For $\la\in\Bbb C$, let $F_{\la}$ denote the
$\bold{a}$-module defined by one generator $v_{\la}$ and relations
$\bold{a}_+v_{\la}=0$, $b_0v_{\la}=2\la v_{\la}$,
$\bold{1} v_{\la}=v_{\la}$.

\sbh{8.5.}Let us introduce the generating functions
$$
\gather
\phi(z)=q-b_0\log\(z)+\sum_{n\ne 0}\frac{b_n}{n}\,z^{-n},\qquad
p(z)=\phi'(z)=-\sum_nb_nz^{-n-1},\\
\beta(z)=\sum_na_nz^{-n-1},\qquad\ga(z)=\sum_n a^*_nz^{-n}\.
\endgather
$$
We have
$$
\alignat2
\phi(z\)\phi(w)&=2\log\(z-w)+\cdots,&\qquad
p(z\)\phi(w)&=\frac{2}{z-w}+\cdots,\\
p(z\)p(w)&=\frac{2}{(z-w)^2}+\cdots,&\qquad
\ga(z\)\beta(w)&=\frac{1}{z-w}+\cdots,
\endalignat
$$
all other products being trivial (do not have a singular part).

\sbh{8.6.}Let us define the currents
$$
\align
E(z)&=\beta(z),\tag"(a)"\\
H(z)&=2\,\Y{\ga(z\)\beta(z)}+\nu p(z),\tag"(b)"\\
F(z)&=-\,\Y{\ga(z)^2\beta(z)}-\nu\,\Y{\ga(z\)p(z)}-k\ga'(z),\tag"(c)"
\endalign
$$
where $\nu^2=k+2$.

\proclaim{\rom{8.7.}\quad Theorem}
The previous formulas define the bosonization
for $\hat{\frak{g}}$, i.e., the Fourier components
of the currents $E(z)$, $H(z)$, $F(z)$
satisfy the commutation relations of $\hat{\frak{g}}$.
\endproclaim

\demo{Proof} We must check the relations (a)--(f) below.
$$
H(z\)H(w)=\frac{2k}{(z-w)^2}+\cdots\.\tag"(a)"
$$
Indeed, we have (using Wick's theorem)
$$
\align
H(z\)H(w)&=(2\,\Y{\ga(z\)\beta(z)}+\nu p(z))(2\,\Y{\ga(w\)\beta(w)}
+\nu p(w))\\
&=4\,\Y{\ga(z\)\beta(z)}\,\Y{\ga(w\)\beta(w)}+\nu^2p(z\)p(w)+\cdots\\
\intertext{(the terms of the first order cancel out)}
&=4\>\{\ga(z\)\beta(w)\}\{\beta(z\)\ga(w)\}+\nu^2p(z\)p(w)+\cdots\\
&=\frac{-4+2\nu^2}{(z-w)^2}+\cdots=\frac{2k}{(z-w)^2}+\cdots\.
\endalign
$$
$$
H(z\)E(w)=\frac{2E(w)}{z-w}+\cdots\.\tag"(b)"
$$
Indeed,
$$
H(z\)E(w)=(2\,\Y{\ga(z\)\beta(z)}
+\nu p(z)\)\beta(w)=2\>\{\ga(z\)\beta(w)\}\>\beta(w)+
\cdots=\frac{2\beta(w)}{z-w}+\cdots\.
$$
$$
H(z\)F(w)=-\frac{2F(w)}{(z-w)}+\cdots\.\tag"(c)"
$$
Indeed,
$$
H(z\)F(w)=-(2\,\Y{\ga(z\)\beta(z)}+\nu p(z))(\,\Y{\ga(w)^2\beta(w)}
+\nu\,\Y{\ga(w\)p(w)}+k\ga'(w))\.
$$
Let us compute all the nontrivial products. We have
$$
\split
\,\Y{\ga(z\)\beta(z)}\,\Y{\ga(w)^2\beta(w)}
&=-\frac{2}{(z-w)^2}\,\ga(w)-\frac{1}{z-w}\,\Y{\beta(w)\ga(w)^2}+\cdots,\\
\,\Y{\ga(z\)\beta(z)}\,\Y{\ga(w\)p(w)}
&=-\frac{1}{z-w}\,\,\Y{\ga(w\)p(w)}+\cdots,\\
\,\Y{\ga(z\)\beta(z)}\ga'(z)&=-\frac{1}{(z-w)^2}\,\ga(z)+\cdots\\
&=-\frac{1}{(z-w)^2}\,\ga(w)-\frac{1}{z-w}\,\ga'(w)+\cdots,\\
p(z)\,\Y{\ga(w\)p(w)}&=\frac{2}{(z-w)^2}\,\ga(w)+\cdots\.
\endsplit
$$
Summing, we get (c).
$$
E(z\)E(w)=0+\cdots\.\tag"(d)"
$$
This is obvious.
$$
E(z\)F(w)=\frac{k}{(z-w)^2}+\frac{H(w)}{z-w}+\cdots\.\tag"(e)"
$$
Indeed,
$$
\split
E(z\)F(w)&=-\beta(z)(\,\Y{\ga(w)^2\beta(w)}
+\nu\,\Y{\ga(w\)p(w)}+k\ga'(w))\\
&=\frac{2}{z-w}\,\,\Y{\ga(w\)\beta(w)}
+\frac{\nu}{z-w}\,p(w)+\frac{k}{(z-w)^2}+\cdots\\
&=\frac{k}{(z-w)^2}+\frac{H(w)}{z-w}+\cdots\.
\endsplit
$$
$$
F(z\)F(w)=0+\cdots\.\tag"(f)"
$$
Indeed,
$$
\split
F(z\)F(w)&=(\,\Y{\ga(z)^2\beta(z)}+\nu.\,\Y{\ga(z\)p(z)}+k\ga'(z))\\
&\qquad\times(\,\Y{\ga(w)^2\beta(w)}+\nu\,\Y{\ga(w\)p(w)}+k\ga'(w))\.
\endsplit
$$
Let us compute the nontrivial products.
$$
\allowdisplaybreaks
\align
\,\Y{\ga(z)^2\beta(z)}\,\Y{\ga(w)^2\beta(w)}
&=-\frac{4}{(z-w)^2}\,\,\Y{\ga(z\)\ga(w)}
+\frac{2}{z-w}\,\,\Y{\ga(z\)\beta(z\)\ga(w)^2}\\
&\qquad-\frac{2}{z-w}\,\,\Y{\ga(z)^2\ga(w\)\beta(w)}+\cdots\\
&=-\frac{4}{(z-w)^2}\,\,\Y{\ga(w)^2}-\frac{4}{z-w}\,\,\Y{\ga(w\)\ga'(w)}
+\cdots,\\
\,\Y{\ga(z\)p(z)}\,\Y{\ga(w\)p(w)}
&=\frac{2}{(z-w)^2}\,\,\Y{\ga(z\)\ga(w)}+\cdots\\
&=\frac{2}{(z-w)^2}\,\,\Y{\ga(w)^2}+\frac{2}{z-w}\,\,\Y{\ga(w\)\ga'(w)}+
\cdots,\\
\,\Y{\ga(z)^2\beta(z)}\,\Y{\ga(w\)p(w)}
&=-\frac{1}{z-w}\,\,\Y{\ga(w)^2p(w)}+\cdots,\\
\,\Y{\ga(z\)p(z)}\,\Y{\ga(w)^2\beta(w)}&=\frac{1}{z-w}\,\,\Y{p(w\)\ga(w)^2}
+\cdots,\\
\,\Y{\ga(z)^2\beta(z)}\ga'(w)&=-\frac{1}{(z-w)^2}\,\,\Y{\ga(z)^2}+\cdots\\
&=-\frac{1}{(z-w)^2}\,\,\Y{\ga(w)^2}-\frac{2}{z-w}\,\,\Y{\ga'(w\)\ga(w)}+\cdots,
\\
\ga'(z)\,\Y{\ga(w)^2\beta(w)}&=-\frac{1}{(z-w)^2}\,\,\Y{\ga(w)^2}+\cdots\.
\endalign
$$
After summing, we get a zero singular part. This proves (f)
and completes the proof of the theorem.\qed
\enddemo

\head\S9.\enspace  Screening current\endhead

\sbh{9.1.}We keep the assumptions of the previous section,
in particular $\frak{g}=\frak{sl}(2)$.
We will assume throughout that $\nu\ne 0$ (the level is noncritical).
The results of this section are due essentially to Feigin-Frenkel
\cite{FeFr1}, \cite{FeFr2}

For $\chi\in\Bbb C$, we shall denote by $W_{\chi;\nu}$
the $\hat{\frak{g}}$-module
$F_{-\chi/2\nu}$, the $\hat{\frak{g}}$-module structure being defined
by formulas 8.6 (a)--(c). Such $\hat{\frak{g}}$-modules are
called {\it Wakimoto modules}.

\sbh{9.2.}For $\a\in\Bbb C$, define the operator
$$
V(\a;z)=\,\Y{\!\exp\(-\a\phi(z))}=T_\a z^{\a b_0}
\exp\bigg(-\a\sum_{n<0}\frac{b_n}{n}\,z^{-n}\bigg)
\exp\bigg(-\a\sum_{n>0}\frac{b_n}{n}z^{-n}\bigg)
$$
acting from $F_{\la}$ to $F_{\la+\a}\otimes z^{2\a\la}\Bbb C((z^{-1}))$.

\proclaim{\rom{9.3.}\quad Lemma} We have
$$
\split
F(z)(-\beta(w\)V(\a;w))
&=-\frac{2\a\nu+k}{(z-w)^2}\,V(\a;w)+
\frac{2-2\a\nu}{z-w}\,\,\Y{\ga(w\)\beta(w\){:}V(\a;w)}\\
&\qquad+\frac{\nu}{z-w}\,\,\Y{p(w\)V(\a;w)}+\cdots\.
\endsplit
$$
\endproclaim

\demo{Proof} The right hand side is equal to
$$
(\,\Y{\ga(z)^2\beta(z)}+\nu\,\,\Y{\ga(z\)p(z)}+k\ga'(z))(\beta(w\)V(\a;w))\.
$$
Let us compute all the products using the Wick theorem. We have
$$
\,\Y{\ga(z)^2\beta(z)}\beta(w\)V(\a;w)=\frac{2}{z-w}\,\,\Y{\ga(w\)\beta(w)}
V(\a;w)+\cdots;
$$
recalling that by 6.4 (a)
$$
p(z\)V(\a;w)=-\frac{2\a}{z-w}\,V(\a;w)+\cdots,
$$
we have
$$
\split
\,\Y{\ga(z\)p(z)}\beta(w\)V(\a;w)
&=-\frac{2\a}{(z-w)^2}\,V(\a;w)+\frac{1}{z-w}\,\,\Y{p(w\)V(\a;w)}\\
&\qquad-\frac{2\a}{z-w}\,\,\Y{\ga(w\)\beta(w)}V(\a;w)+\cdots\.
\endsplit
$$
Finally,
$$
\ga'(z\)\beta(w\)V(\a;w)=-\frac{1}{(z-w)^2}\,V(\a;w)+\cdots\.
$$
By summing all up, we get the statement of the lemma.\qed
\enddemo

\sbh{9.4.}Let us introduce the operators called {\it screening currents}
by
$$
S(z)=-\beta(z\)V(\nu^{-1};z)\:\enspace W_{\chi;\nu}\to W_{\chi-2;\nu}\otimes
z^{-\chi/\nu^2}\Bbb C((z^{-1}))\qquad (\chi\in\Bbb C)\.
$$
Set
$$
S(F;z)=-\nu^2V(\nu^{-1};z)\.
$$

\proclaim{\rom{9.5.}\quad Theorem} We have
$$
F(z\)S(w)=\paw \bigg(\frac{S(F;w)}{z-w}\bigg)+\cdots,\quad
E(z\)S(w)=0+\cdots,\quad H(z\)S(w)=0+\cdots\.
$$
\endproclaim

\demo{Proof} The first formula follows from the previous lemma after
the substitution $\a=\nu^{-1}$, taking into account that
$$
\paw V(\a;w)=-\a\,\Y{p(w\)V(\a;w)}\.
$$
The second formula is obvious. Let us prove the third formula.
We have
$$
H(z\)S(w)=-(2\,\Y{\ga(z\)\beta(z)}+\nu p(z)\)\beta(w\)V(\nu^{-1};w)\.
$$
Now,
$$
\align
\ga(z\)\beta(z\)\beta(w)V(\nu^{-1};w)
&=\frac{1}{z-w}\,\beta(w\)V(\nu^{-1};w)+\cdots,\\
p(z\)\beta(w\)V(\nu^{-1};w)&=-\frac{2\nu^{-1}}{z-w}\,
\beta(w\)V(\nu^{-1};w)+\cdots\.
\endalign
$$
By adding the two terms, we get the third formula.\qed
\enddemo

\proclaim{\rom{9.6.}\quad Corollary} For every $n\in\Bbb Z$,
$$
[F_n,S(w)]=\paw (w^nS(F;w)),\qquad[E_n,S(w)]=[H_n,S(w)]=0\.\qed
$$
\endproclaim

\sbh{9.7.}We define the operators $S(X;z)$ (\<$X\in\frak{g}$) as
follows. First, $S(X;z)$ linearly depends on $X\in\frak{g}$. Second,
$S(F;z)$ is as above, and we put $S(E;z)=S(H;z)=0$.

Now,  we define the operators $S(x;z)$ (\<$x\in\hat{\frak{g}}$) as
follows. First,
they depend linearly on $x\in\hat{\frak{g}}$. Second, we set
$$
S(X_n;z)=z^nS(X;z)\quad(X\in\frak{g},\;n\in\Bbb Z),\qquad
S(\bold{1};z)=0\.\tag"(a)"
$$
In this notation, we can rewrite the previous theorem and
corollary as follows.

\proclaim{\rom{9.8.}\quad Theorem}
\rom{(a)} For all $X\in\frak{g}$,
$$
X(z\)S(w)=\paw \bigg(\frac{S(X;w)}{z-w}\bigg)+\cdots\.
$$

\noindent\rom{(b)} For all $x\in\hat{\frak{g}}$,
$$
[x,S(w)]=\paw  S(x;w).\qed
$$
\endproclaim

\proclaim{\rom{9.9.}\quad Lemma} We have
$$
\gather
E(z\)S(F;w)=0+\cdots,\\
H(z\)S(F;w)=-2\,\frac{S(F;w)}{z-w}+\cdots,\qquad
F(z\)S(F;w)=\frac{\ga(w\)S(F;w)}{z-w}+\cdots\.
\endgather
$$
\endproclaim

This is proved by a simple direct computation.

\proclaim{\rom{9.10.}\quad Theorem}
\rom{(a)} For any $X,Y\in\frak{g}$, we have
$$
X(z\)S(Y;w)-Y(z\)S(X;w)=\frac{S([X,Y];w)}{z-w}+\dots\.
$$

\noindent\rom{(b)} For any $x,y\in\hat{\frak{g}}$, we have
$$
[x,S(y;w)]-[y,S(x;w)]=S([x,y];w)\.
$$
\endproclaim

\demo{Proof} (a) follows from the previous lemma. Alternatively, it follows
from 9.8 (a) and the facts (c), (d) below, using the associativity
of the operator products.

(b) follows from (a) and the fact that

\itm{c}{in the products $X(z\)S(Y;w)$
at most first order poles are present.}

(This claim is a weaker version of the previous lemma.)

Alternatively,
(b) follows from 9.8 (b), if we notice that

\itm{d}{for generic $\chi$ the operator $\paw $ is an isomorphism.}

This proves the theorem.\qed
\enddemo

\sbh{9.11.}Set $A=\Bbb C((z^{-1}))$. Consider the twisted
de Rham complex
$$
\Om^{\bdot}\:\enspace 0@>>>\Om^0@>>>\Om^1@>>>0,
$$
where $\Om^0=A$, $\Om^1=A\,dz$, the differential being equal to
$$
d_{DR}-\frac{\chi}{\nu^2}\cdot\frac{dz}{z}\.
$$
As before, we form the double complex
$$
C^{\bdot}(\hat{\frak{g}};\op{Hom}\(W_{\chi;\nu},W_{\chi-2;\nu}\otimes
\Om^{\bdot}))\.
$$
Let us define a one-cochain $V=(V^{01},V^{10})$ in the associated
total complex as follows. We set
$$
V^{01}=S(z\)z^{\chi/\nu^2}dz,\quad
V^{10}(x)=S(x;z\)z^{\chi/\nu^2}\qquad (x\in\hat{\frak{g}})\.
$$

The next theorem is the reformulation of Theorems 9.8 (b) and
9.10 (b).

\proclaim{\rom{9.12.}\quad Theorem} The cochain $V$ is a $1$-cocycle.\qed
\endproclaim

\sbh{9.13.}For an integer $p\ge 1$, consider the normally ordered product
$$
\,\Y{S(z_1)\cdots S(z_p)}\.
$$
To simplify the notations, we regard it as an element of
$$
\op{Hom}\(W_{\chi;\nu},W_{\chi-2p;\nu}\otimes A_p)
$$
(we ignore the powers of $z_i$). The same concerns the operators
$S(x;z_i)$.
Recall that
$$
A_p=\Bbb C[[z_1,\dots,z_p]]
\bigg[\prod z_i^{-1}, \prod_{i\ne j}(z_i-z_j)^{-1}\bigg]=\Om^0(X_p)
$$
(cf. 7.4).

Consider the twisted de Rham complex
$$
\Om^{\bdot}\:\enspace 0@>>>\Om^0@>>>\cdots@>>>\Om^p@>>>0,
$$
where $\Om^i=\Om^i(X_p)$. The differential is
$$
d_{DR}-\sum_i\frac{\chi}{\nu^2}\cdot \frac{dz_i}{z_i}+
\sum_{i<j}\frac{2}{\nu^2}\cdot\frac{dz_i-dz_j}{z_i-z_j}\.
$$

\sbh{9.14.}For each $a=0,\dots,p$,  we define the operators
$$
V^{a,p-a}\in\op{Hom}\(\Lambda^a\hat{\frak{g}};\op{Hom}\(W_{\chi;\nu},W_{\chi
-2p;\nu}\otimes
\Om^{p-a}))
$$
as follows. By definition,
$$
\split
&V^{a,p-a}(x_1\wedge\cdots\wedge x_a)\\
&\qquad=\sum_{1\le i_1<\cdots<i_a\le p}(-1)^{\op{sgn}(i_1,\dots,i_a)}\\
&\qquad\qquad\times\bigg(\sum_{\sigma\in\Sigma_m}(-1)^{\op{sgn}(\sigma)}
\,\Y{S(z_1)\cdots
S(x_{\sigma(1)};z_{i_1})\cdots S(x_{\sigma(a)};z_{i_a})\cdot
\dots\cdot S(z_p)}\bigg)\\
&\qquad\qquad\times dz_1\wedge\cdots\wedge\hat{dz}_{i_1}\wedge\cdots
\wedge\hat{dz}_{i_a}\wedge\cdots\wedge dz_p\.
\endsplit
$$
Here the sign $\op{sgn}\(i_1,\dots,i_a)$ is
defined by induction on $a$ as follows.
$$
\op{sgn}\(\;\;)=0,\qquad\op{sgn}(i_1,\dots,i_a)
=\op{sgn}(i_1,\dots,i_{a-1})+i_a+a\.
$$
For example,
$$
\align
V^{0p}&=\,\Y{S(z_1)\cdots S(z_p)}\,dz_1\wedge\cdots\wedge dz_p,\\
V^{p0}(x_1\wedge\cdots\wedge x_p)
&=\sum_{\sigma\in\Sigma_p}(-1)^{\op{sgn}\(\sigma)}
\,\Y{S(x_{\sigma(1)};z_1)\cdots S(x_{\sigma(p)};z_p)}\.
\endalign
$$

\sbh{9.15.}Consider the double Koszul complex $C^{\bdot}
(\hat{\frak{g}};\op{Hom}\(W_{\chi;\nu},W_{\chi-2p;\nu}\otimes\Om^{\bdot}))$.

We have a $p$-cochain $V=(V^{0p},\dots,V^{p0})$ in the associated
total complex.

\proclaim{\rom{9.16.}\quad Theorem} The cochain $V$ is a $p$-cocycle.\qed
\endproclaim

\head Chapter 4. Affine Lie algebras (general case)\endhead

\head\S10.\enspace  Bosonization\endhead

The basic mathematical results on bosonization for
arbitrary affine Lie algebras are due to B. Feigin and E. Frenkel,
see~\cite{FeFr1}, \cite{FeFr2}, \cite{Fr}, and also
\cite{FFR} and references therein. For an original
physical approach see \cite{BMP2}.

\sbh{10.1.} 
Let $\frak{g}$ be the finite-dimensional  complex
Lie algebra\footnote{We suspect that most of what appears below is true for
the Kac--Moody algebra
corresponding to an arbitrary symmetrizable generalized Cartan matrix.}
corresponding to a Cartan matrix $A=(a_{ij})_{i,j=1}^r$ and
with the Chevalley generators
$E_i$, $H_i$, $F_i$ (\<$i=1,\dots, r$).
Let $\frak{h}\subset\frak{g}$ be the Cartan subalgebra
generated by $H_1,\dots, H_r$.
Let $\a_1,\dots,\a_r\in\frak{h}^*$ be the simple roots; let
$(\bcdot,\bcdot)$ be the symmetric nondegenerate bilinear form on
$\frak{h}^*$
such that $a_{ij}=2\(\a_i,\a_j)/(\a_i,\a_i)$.
This bilinear form defines an isomorphism $\frak{h}^*@>\sim>>\frak{h}$;
using this
isomorphism, we carry over the bilinear form to $\frak{h}$;
this last bilinear
form will also be denoted by $(\bcdot,\bcdot)$. Finally, $\Delta_+$
will denote the set of positive roots.

Let $g$ be the dual Coxeter number of the root system of $\frak{g}$. We fix
a
complex parameter $\nu\ne 0$ and set $k=\nu^2-g$.

\sbh{10.2.}Let $\bold{a}$ be the Lie algebra defined by the generators
$b_n^i$ (\<$i=1,\dots, r$, $n\in\Bbb Z$), $a_n^{\a}$, $a_n^{\a*}$
(\<$\a\in\Delta_+$, $n\in\Bbb Z$), $\bold{1}$ and the relations

\itm{a}{$[b^i_n,b^j_m]=(H_i,H_j\)n\delta_{ij}\delta_{n+m,0}\cdot\bold{1}$;}

\itm{b}{$[a^{\a}_n,a^{\beta*}_m]=\delta_{\a\beta}\delta_{n+m,0}\cdot\bold{1}
$;}

\itm{c}{all other commutators between the generators vanish.}

Let $\bold{a}_+$ denote the Lie subalgebra of $\bold{a}$ generated
by the elements
$b_n^i$, $a_n^{\a}$ (\<$n>0$, $\a\in\Delta_+$), $a^{\a*}_n$
(\<$n\ge 0$, $\a\in\Delta_+$). These generators are called
{\it annihilation operators}. One introduces the {\it normal ordering\/}
of a monomial in $\U\bold{a}$
in the usual way: all annihilation operators should be pulled to the right.

All the $\bold{a}$-modules $M$ that we shall
consider will be {\it $\bold{a}_+$-locally finite}, i.e., have the property
that for every $x\in M$, the space $\U\bold{a}_+x$ is finite-dimensional.
The operator $\bold{1}$ will act as the identity.

For computational purposes, we shall also use the operators $q^i$
(\<$i=1,\dots,r$),
with the only nontrivial commutation relations
$$
[q^i,b^j_0]=(H_i,H_j)\cdot\bold{1}\.
$$

\sbh{10.3.}For $\la\in\frak{h}^*$, let $\Cal{F}_{\la}$ denote the
$\bold{a}$-module defined by one generator $v_{\la}$ and the relations
$\bold{a}_+v_{\la}=0$, $b_0^iv_{\la}=\langle H_i,\la\rangle\,
\la v_{\la}$, $\bold{1} v_{\la}=v_{\la}$.

\sbh{10.4.}Let us introduce the generating functions
$$
\gather
\phi^i(z)=q^i-b^i_0\log\(z)+\sum_{n\ne 0}\frac{b^i_n}{n}\,z^{-n},\\
p^i(z)=\phi^{i\prime}(z)=-\sum_nb^i_nz^{-n-1}\qquad (i=1,\dots,r),\\
\beta^{\a}(z)=\sum_na^{\a}_nz^{-n-1},\quad
\ga^\a(z)=\sum_na^{\a*}_nz^{-n}\qquad(\a\in\Delta_+)\.
\endgather
$$
We have
$$
\phi^i(z\)\phi^j(w)=(H_i,H_j)\log\(z-w)+\cdots,\qquad
\ga^{\a}(z\)\beta^{\a'}(w)=\frac{\delta_{\a\a'}}{z-w}+\cdots,
$$
all other products being trivial.

To shorten the notation, we shall write $\beta^i(z)$,
$\ga^i(z)$ instead of $\beta^{\a_i}(z)$, $\ga^{\a_i}(z)$.

\sbh{10.5.}The bosonization formulas for the affine Lie algebra
$\hat{\frak{g}}$ with central charge $k$ have the following form.
$$
\align
E_i(z)&=\,\Y{\Cal{E}_i(\ga^{\a}(z),\beta^{\a'}(z))},\tag"(a)"\\
H_i(z)&=\sum_{\a\in\Delta_+}\langle H_i,\a\rangle
\,\Y{\ga^{\a}(z\)\beta^{\a}(z)}+\nu p^i(z),\tag"(b)"\\
F_i(z)&=\,\Y{\Cal{F}_i(\ga^{\a}(z),\beta^{\a'}(z))}-
\nu\,\frac{(\a_i,\a_i)}{2}\,\Y{\ga^i(z\)p^i(z)}-c_i(\nu)\,
\frac{d\ga^i(z)}{dz}\.\tag"(c)"
\endalign
$$
Here $\Cal{E}_i(\ga^{\a},\beta^{\a'})$, $\Cal{F}_i(\ga^{\a},
\beta^{\a'})$ are certain polynomials, linear in the $\beta$'s.
They depend on the Poincar\'{e}--Birkhoff--Witt isomorphism
for the enveloping algebra $\U\frak{n}_-$, which identifies this
algebra with
the symmetric algebra on root vectors $F_\a$ (\<$\a\in\Delta_+$).
For their definition, see \cite{BMP2},\cite{FeFr1}, \cite{FeFr2}.
 The coefficients
$c_i(\nu)$ are certain numbers.

\head\S11.\enspace  Screening currents\endhead

\sbh{11.1.}For $\chi\in\frak{h}^*$, we define the Wakimoto
module $W_{\chi;\nu}$ as the $\hat{\frak{g}}$-module
$\Cal{F}_{-\chi/\nu}$, the
$\hat{\frak{g}}$-module structure being defined by the bosonization
formulas
10.5 (a)--(c).

\sbh{11.2.}For $\mu=\sum_i\mu_i\a_i\in\frak{h}^*$, define the
bosonic vertex operator
$$
V(\mu;z)=\,\Y{\!\exp\bigg(-\sum_i\mu_i\,\frac{(\a_i,\a_i)}{2}\,\phi^i(z)\bigg)
}\.
$$
This operator acts from $\Cal{F}_{\la}$ to $\Cal{F}_{\la+\mu}\otimes
z^{(\la,\mu)}\Bbb C((z^{-1}))$.

\sbh{11.3.}We have the {\it product formula\/} (cf. 6.11)
$$
V(\mu_1;z_1\)V(\mu_2;z_2)=(z_1-z_2)^{(\mu_1,\mu_2)}
\,\Y{V(\mu_1;z_1\)V(\mu_2;z_2)}\.
$$

\sbh{11.4.}By definition, the {\it screening currents\/} are defined
by
$$
S_i(z)=\Cal{S}_i(\ga^{\a}(z),\beta^{\a'}(z)\)V(\nu^{-1}\a_i;z)\:\enspace 
W_{\chi;\nu}\to W_{\chi-\a_i;\nu}\otimes z^{-(\chi,\a_i)/\nu^2}
\Bbb C((z^{-1}))\.
$$
Here $\Cal{S}_i(\ga^{\a},\beta^{\a'})$ are certain polynomials
depending on a PBW decomposition for $\Cal{U}\frak{n}_-$, see
\cite{BMP2}.

We set
$$
S_i(F_j;z)=-\delta_{ij}\nu^2V(\nu^{-1}\a_i;z)\.
$$

\sbh{11.5.} The following operator expansion formulas which were
proved in [BMP1], [FeFr1] and [FeFr2] for all classical simple Lie
algebras, essentially on a case by case basis, will play a key role below
$$
\gather
F_i(z\)S_j(w)=\paw \bigg(\frac{S_j(F_i;w)}{z-w}\bigg)+\cdots,\\
E_i(z\)S_j(w)=0+\cdots,\qquad H_i(z\)S_j(w)=0+\cdots\.
\endgather
$$

\sbh{11.6.}Let ${\bold g}$ be the free Lie algebra with generators
$E_i,\, F_i$ and $H_i,\,i=1,\ldots,r.$ To each $i=1,\ldots,r$ and
$X\in {\bold g}$ we associate an operator $S_i(X;z)$ as follows.
On generators we put $S_i(E_j;z)=S_i(H_j;z)=0$, and let
$S_i(F_j;z)$ be as above. Then define the $S_i(X;z)$
inductively by the formula
$$
[X,S_i(Y;w)]-[Y,S_i(X;w)]=S_i([X,Y];w)\quad,\quad \forall X,Y\in {\bold g}\,
$$
and by $\Bbb C$-linearity.
We further extend this definition to the loop Lie algebra
${\bold g}((z))$ by setting
$$
S_i(X_n;z)=z^nS(X;z)\quad,\quad
i=1,\dots,r\;,\,X\in{\bold g},\;n\in\Bbb Z\.\tag"(a)"
$$

\sbh{11.7.}Let ${\bold g} \twoheadrightarrow \frak{g}$ be the canonical
projection from the free Lie algebra to the semisimple
Lie algebra. Write $\pi: {\bold g}((z)) \to \hat{\frak{g}}$
for the induced map of on loops. It is straightforward to deduce
by induction from formulas of n.11.5, that 

 {\it For all $x\in{\bold g}((z))$, $i=1,\dots, r$, we have}
$$
[\pi(x),S_i(w)]=\paw  S_i(x;w)\.\tag 11.7.1
$$

\proclaim{\rom{11.8.}\quad Theorem} The assignment $x\mapsto S_i(x;w)$
descends, for any $i=1,\dots, r$, to a well-defined map
$\hat{\frak{g}} \to \op{Hom}\bigl(
W_{\chi;\nu}\,,\, W_{\chi-\a_i;\nu}\otimes z^{-(\chi,\a_i)/\nu^2}
\Bbb C((z^{-1}))\bigr)\.$ The following equations hold for any
$x,y\in\hat{\frak{g}}, \,i=1,\dots, r$
$$
[x,S_i(y;w)]-[y,S_i(x;w)]=S_i([x,y];w)\quad,\quad
[x,S_i(w)]=\paw  S_i(x;w)\.
$$
\endproclaim

\demo{Proof} The only statement that is not immediate from construction
is that the map $x\mapsto S_i(x;w)$ descends from ${\bold g}((z))$
to $\hat{\frak{g}}$. To prove this, let $x,x'\in {\bold g}((z))$ be such
that $\pi(x)=\pi(x') \in  \hat{\frak{g}}$.
Then equation (11.7.1) yields
$$\paw  S_i(x;w) = \paw  S_i(x';w)\.$$
But for generic
$\chi$ the operator $\paw $ is an isomorphism on the de Rham complex
of our local system. Hence, for generic
$\chi$, we can conclude that $S_i(x;w) =  S_i(x';w)$.
Now, for arbitrary $\chi$, the equation $S_i(x;w) =  S_i(x';w)$
follows by continuity.\qed
\enddemo

\proclaim{\rom{11.9.}\quad Conjecture}
For any $X,Y\in\frak{g}$, $i=1,\dots, r$, one has
$$
X(z\)S_i(Y;w)-Y(z\)S_i(X;w)=\frac{S_i([X,Y];w)}{z-w}+\cdots\.
$$
\endproclaim

\sbh{11.10.}Let us fix a sequence $i_1,\dots,i_p$, where
$1\le i_j\le r$ for all $j$.

Consider the twisted de Rham complex
$$
\Om^{\bdot}\:\enspace 0@>>>\Om^0@>>>\cdots@>>>\Om^p@>>>0,
$$
where $\Om^i=\Om^i(X_p)$ are the same spaces as in 9.13.
The differential is by definition equal to
$$
d_{DR}-\sum_{j=1}^p\frac{(\chi,\a_{i_j})}{\nu^2}\cdot \frac{dz_j}{z_j}+
\sum_{1\le j'<j''\le p}\frac{(\a_{i_{j'}},\a_{i_{j''}})}{\nu^2}
\cdot\frac{dz_{j'}-dz_{j''}}{z_{j'}-z_{j''}}\.
$$

\sbh{11.11.}Set $\a=\sum_{j=1}^p\a_{i_j}$.
For each $a=0,\dots,p$,  we define the operators
$$
V^{a,p-a}\in\op{Hom}\(\Lambda^a\hat{\frak{g}};
\op{Hom}\(W_{\chi;\nu},W_{\chi-\a;\nu}\otimes\Om^{p-a}))
$$
as follows. By definition,
$$
\split
&V^{a,p-a}(x_1\wedge\cdots\wedge x_a)\\
&\quad=\sum_{1\le j_1<\cdots<j_a\le p}(-1)^{\op{sgn}(j_1,\dots,j_a)}\\
&\qquad\times\bigg(\sum_{\sigma\in\Sigma_m}(-1)^{\op{sgn}(\sigma)}
\,\Y{S_{i_1}(z_1)\cdots S_{i_{j_1}}(x_{\sigma(1)};z_{j_1})
\cdots S_{i_{j_a}}(x_{\sigma(a)};z_{j_a})\cdots S_{i_p}(z_p)}\bigg)\\
&\qquad\times dz_1\wedge\dots\wedge\hat{dz}_{j_1}\wedge\dots\wedge
\hat{dz}_{j_a}\wedge\dots\wedge dz_p
\endsplit
$$
(cf. 9.14). The sign $\op{sgn}(j_1,\dots,j_a)$ is defined in 9.14.

For example,
$$
\align
V^{0p}&=\,\Y{S_{i_1}(z_1)\cdots S_{i_p}(z_p)}\,
dz_1\wedge\dots\wedge dz_p,\\
V^{p0}(x_1\wedge\cdots\wedge x_p)
&=\sum_{\sigma\in\Sigma_p}(-1)^{\op{sgn}(\sigma)}
\,\Y{S_{i_1}(x_{\sigma(1)};z_1)\cdots S_{i_p}(x_{\sigma(p)};z_p)}\.
\endalign
$$

\sbh{11.12.}Consider the double Koszul complex 

$C^{\bdot}
(\hat{\frak{g}};\op{Hom}\(W_{\chi;\nu},W_{\chi-\a;\nu}\otimes\Om^{\bdot}))$.

We have a $p$-cochain $V=(V^{0p},\dots,V^{p0})$ in the associated
total complex.

\proclaim{\rom{11.13.}\quad Theorem}
The cochain $V$ is a $p$-cocycle.\qed
\endproclaim


\Refs\widestnumber\key{BMP2}

\ref\key{BD} 
\by A. Beilinson, V. Drinfeld
\paper Chiral algebras
\jour Preprint
\yr 1996
\endref

\ref\key{BMP1}
\by P. Bowknegt, J. McCarthy, and K. Pilch
\paper Quantum group structure in the Fock space resolutions of the
$\widehat{sl}(n)$ representations
\jour Comm. Math. Phys.
\vol 131
\yr 1990
\pages 125--155
\endref

\ref\key{BMP2}
\bysame
\paper Free field approach to two-dimensional conformal field theory
\jour Prog. Theor. Phys. Suppl.
\vol 102
\yr 1990
\pages 67--135
\endref

\ref\key{Br}
\by J.-L. Brylinski
\paper Non-commutative Ruelle-Sullivan type currents
\inbook in: The Grothendieck Festschrift Volume I
\eds P. Cartier et al. 
\publ Birkh\"auser
\publaddr Boston et al. 
\yr 1990
\pages 477--498
\endref 

\ref\key{FFR}
\by B. Feigin, E. Frenkel, and N. Reshetikhin
\paper Gaudin model, Bethe Ansatz and critical level
\jour Comm. Math. Phys. 
\vol 166
\yr 1994
\pages 27--62
\endref

\ref\key{FeFr1}
\by B. Feigin and E. Frenkel
\jour Comm. Math. Phys.
\vol 128
\yr 1990
\pages 161--189;
{\it Comm. Math. Phys.} {\bf 137} (1992), 617--639;
{\it  Phys. Letters.} {\bf B246} (1990), 75--81;
{\it   Lett. Mathem. Physics.} {\bf 19} (1990),
307--317
\endref

\ref\key{FeFr2}
\by B. Feigin and E. Frenkel
\inbook in: V. Knizhnik memorial volume
\eds L. Brink, D. Friedan, A. Polyakov
\publ World Scientific
\publaddr Singapore
\yr 1990
\pages 271--316
\endref

\ref\key{FF}
\by B. Feigin and D. Fuchs
\paper Representations of the Virasoro
algebra
\inbook in: Representations of Lie groups and related topics
\eds A. M. Vershik and D. P. Zhelobenko
\publ Gordon and Breach
\publaddr New York
\yr 1990
\pages 465--554
\endref


\ref\key{F}
\by G. Felder
\paper BRST approach to minimal models
\jour Nuclear Phys. B
\vol 317
\yr 1989
\pages 215--236
\endref

\ref\key{FS}
\by M. Finkelberg and V. Schechtman 
\paper Localization of modules over small quantum groups
\jour J. of Math. Sciences 
\vol 82 
\yr 1996
\pages 3127--3164
\endref

\ref\key{Fr}
\by E. Frenkel
\paper Free field realizations in representation
theory and conformal field theory
\jour Address to the ICM, Z\"{u}rich,
August 3--11
\yr 1994
\endref

\ref\key{K}
\by V. Kac
\book Infinite dimensional Lie algebras
\bookinfo 2nd Edition
\publ Cambridge University Press
\publaddr Cambridge et al.
\yr 1985
\endref

\ref\key{Ka}
\by M. Kashiwara
\paper On crystal bases of the $q$-analogue of the universal enveloping algebra
\jour Duke Math. J.
\vol 63
\yr 1991
\pages 465-516
\endref

\ref\key{L}
\by G. Lusztig
\book Introduction to quantum groups
\publ Birkh\"auser
\publaddr Boston et al.
\yr 1993
\endref

\ref\key{S}
\by A. Sebbar
\paper Quantum screening operators and $q$-de Rham cohomology
\jour PhD Thesis, Stony Brook
\yr 1997
\endref

\ref\key{TK}
\by A. Tsuchiya and Y. Kanie
\paper Fock space representations of the Virasoro algebra
\jour Publ. Res. Inst. Math. Sci.
\vol 22
\yr 1986
\pages 259--327
\endref

\ref\key{W}
\by M. Wakimoto
\paper Fock representations of the affine Lie algebra
$A_1^{(1)}$
\jour Comm. Math. Phys.
\vol 104
\yr 1986
\pages 605--609
\endref

\endRefs
\enddocument